\documentclass[12pt]{article}
\usepackage{amsmath}
\usepackage{graphicx}
\usepackage{enumerate}
\usepackage{natbib}
\usepackage{url} 

\newcommand{\blind}{1}

\usepackage{amsmath}
\usepackage{amsthm}
\usepackage{bbm}
\usepackage{bm}
\usepackage{amssymb}
\usepackage{mathrsfs}
\usepackage{dsfont}
\usepackage{float}
\usepackage{booktabs,multirow}
\usepackage{makecell}
\usepackage{subfigure}

\usepackage{enumitem}

\renewcommand{\baselinestretch}{1.5}
\usepackage{natbib}
\makeatletter
\renewcommand\@biblabel[1]{${#1}.$}
\usepackage[colorlinks=true,linkcolor=blue, citecolor=blue, urlcolor=blue]{hyperref}
\hypersetup{hypertex=true,
	colorlinks=true,
	linkcolor=blue,
	anchorcolor=blue,
	citecolor=blue}
\makeatother
\usepackage{color}
\usepackage[font=small,format=plain,labelfont=bf,up,textfont=normal,up,justification=justified,singlelinecheck=false]{caption}

\usepackage{pgfplots}

\usepackage{authblk}
\usepackage{tikz}
\usetikzlibrary{shapes.geometric, arrows}
\DeclareMathOperator*{\argmax}{argmax}

\usepackage[noend]{algpseudocode}
\usepackage{algorithmicx,algorithm}
\usepackage{epstopdf}
\usepackage{multirow}

\newtheorem{theorem}{Theorem}
\newtheorem{corollary}{Corollary}

\newtheorem{proposition}{Proposition}
\newtheorem{assumption}{Assumption}
\newtheorem{remark}{Remark}

\usepackage[normalem]{ulem}

\begin{document}

\def\spacingset#1{\renewcommand{\baselinestretch}%
{#1}\small\normalsize} \spacingset{1}


\if1\blind
{
  \title{\bf Staleness Factors and Volatility Estimation at High Frequencies}
  \author{Xinbing Kong\\
	Southeast University, Nanjing 211189, China\\
	Bin Wu\thanks{Co-first and Corresponding author. Email: bin.w@ustc.edu.cn.\\
	Authors are listed alphabetically.}\\
	University of Science and Technology of China, Hefei 230026, China\\
	Wuyi Ye\\
	University of Science and Technology of China, Hefei 230026, China}
  \maketitle
} \fi

\if0\blind
{
  \bigskip
  \bigskip
  \bigskip
  \begin{center}
    {\LARGE\bf Staleness Factors and Volatility Estimation at High Frequencies}
\end{center}
  \medskip
} \fi

\bigskip
\begin{abstract}
In this paper, we propose a price staleness factor model that accounts for pervasive market friction across assets and incorporates relevant covariates. Using large-panel high-frequency data, we derive the maximum likelihood estimators of the regression coefficients, the nonstationary factors, and their loading parameters. These estimators recover the time-varying price staleness probabilities. We develop asymptotic theory in which both the dimension $d$ and the sampling frequency $n$ tend to infinity. Using a local principal component analysis (LPCA) approach, we find that the efficient price co-volatilities (systematic and idiosyncratic) are biased downward due to the presence of staleness. We provide bias-corrected estimators for both the spot and integrated systematic and idiosyncratic co-volatilities, and prove that these estimators are robust to data staleness. Interestingly, besides their dependence on the dimensionality $d$, the integrated plug-in estimates converge at a rate of $n^{-1/2}$, whereas the LPCA estimates converge at a slower rate of $n^{-1/4}$. This validates the aggregation efficiency achieved through nonlinear, nonstationary factor analysis via maximum likelihood estimation. Numerical experiments justify our theoretical findings. Empirically, we demonstrate that the staleness factor provides unique explanatory power for cross-sectional risk premia, and that the staleness correction reduces out-of-sample portfolio risk.
\end{abstract}

\noindent%
\vfill

\newpage
\spacingset{1.6} 
\section{Introduction}\label{sect:introduction}

Price staleness refers to the phenomenon where asset prices are updated less frequently than expected. Price staleness is commonly attributed to market frictions that impede the continuous incorporation of information into transaction prices. Under no-arbitrage conditions, asset prices typically evolve as semimartingales, exhibiting stochastic continuity in their paths. When the semimartingale is continuously driven by Brownian motions, high-frequency returns scale with the square root of the time lag. However, \citet{bandi2017excess} shows that a large proportion of high-frequency returns are abnormally small (smaller than what continuous semimartingale models imply).

Staleness probability, defined statistically as the relative frequency of zero returns (named ``zeros''), is influenced by two primary factors: low trading volumes and price discretization (\citealt{bandi2020zeros}). This concept provides valuable insights into market frictions and their underlying determinants (particularly liquidity factors). Since \cite{bandi2017excess} first pioneered zero-return analysis using intraday data in continuous-time frameworks, the staleness literature has expanded significantly (c.f., \citealt{bandi2020zeros}; \citealt{kolokolov2020statistical}; \citealt{bandi2024systematic}; \citealt{liu2024bias}; \citealt{zhu2024bivariate}). For ease of presentation, let $t_j$ and $t_{j-1}$ denote two adjacent sampling times. A widely adopted model in financial econometrics specifies the observed log price $\widetilde{Y}_{t_j}$ at time $t_j$ as:
\begin{equation}\label{unimodel}
	\widetilde{Y}_{t_j}=Y_{t_j}(1-B_{t_j})+\widetilde{Y}_{t_{j-1}}B_{t_j}
\end{equation}
with initial value $\widetilde{Y}_{t_0}=Y_{t_0}$, where $B_{t_j}$ is a Bernoulli random variable indicating whether prices update ($B_{t_j}=0$) or remain unchanged ($B_{t_j}=1$). The sluggish price component $\widetilde{Y}_{t_{j-1}}B_{t_j}$ captures price staleness, while $Y_{t}$ is an efficient price semimartingale.

Existing research has primarily focused on univariate series or fixed-dimension multivariate processes. However, \cite{bandi2024systematic} demonstrates systematic components in price-updating delays, revealing cross-sectionally correlated staleness patterns across assets. Consequently, modeling joint staleness probabilities in large asset pools becomes crucial for statistical theory and financial applications. Though the model (\ref{unimodel}) and the large-dimensional extension (\ref{eq:observable price processes}) below were initially developed within the financial literature, their theoretical framework extends naturally to other contexts, such as streaming-data applications with information delays or data-cleaning procedures in which missing observations are imputed by carrying forward the most recent available value until a new update arrives.

Two fundamental questions naturally arise in practical applications. First, to what extent do staleness factors account for the substantial cross-sectional variation observed in high-frequency data? In the context of large-scale asset pricing, assessing the performance of staleness factors as proxies for liquidity is of considerable importance. Second, does data staleness introduce estimation bias in large volatility matrices? In portfolio allocation, inaccurate volatility matrix estimates can amplify out-of-sample risk in mean-variance optimization strategies. These observations motivate our study.

To the best of our knowledge, no existing study has directly addressed the modeling of price staleness in a high-dimensional setting using a large panel of high-frequency data. One notable exception is the work of \cite{bandi2024systematic}, which provides an initial investigation into the existence of price co-staleness and proposes statistical indicators to measure and explain observed empirical patterns. However, that study relies on the restrictive assumption that zero (or near-zero) returns occur simultaneously across all assets at each time stamp. In practice, however, delays in the transmission of liquidity information across assets can occur. While the probability of stale prices for all assets at any given time is positive, simultaneous zeros across all assets are rare, particularly at high frequencies for high-dimensional price processes. Moreover, \cite{bandi2024systematic} assumes that systematic staleness is constant and driven by a single factor. Our empirical analysis reveals that staleness factor series exhibit clear time variation and nonstationary patterns.

In this paper, we formally introduce a novel nonlinear continuous-time model for high-dimensional staleness processes, termed the staleness factor model (SFM). The model specifies staleness probabilities through exogenous covariates and unobservable common factors via a general link function (e.g., logit or probit), offering several key advantages over existing frameworks. First, by modeling staleness probabilities as a function of these covariates and factors, the SFM naturally accounts for \emph{price staleness pervasiveness}. Even when flat prices are not simultaneously observed across all assets, the staleness probability remains positive, making delayed flat-price arrivals interpretable. Second, allowing both the staleness factors and the covariate processes to vary over time makes the model more flexible and better supported by empirical data. Another key difference from existing continuous-time factor models (such as \citealt{ait2017using}; \citealt{pelger2019large}; \citealt{kong2017number,kong2018systematic}) is that, in our model, the price staleness probability process cannot be differenced, since the price staleness probability (the probability that $B_{t_j}=1$) is unobservable. This poses the challenge for inference, because high-frequency global principal component analysis (GPCA) and local principal component analysis (LPCA) methods (see \citealt{kong2023discrepancy}) that rely on differenced semimartingales become inapplicable. We address this challenge to estimate this nonlinear, nonstationary staleness factor model by employing maximum likelihood estimation (MLE). We show that the estimator of the staleness probability has an error bound of the order $(\min(\sqrt{n},\sqrt{d}))^{-1}$. Furthermore, under suitable regularity conditions, the integrated version of the estimator achieves the $n^{-1/2}$ rate, consistent with the efficiency rate of estimated volatility functionals as theoretically underpinned by \cite{jacod2013quarticity}. Notably, the MLE estimator is not subject to biases due to nonlinearity, volatility-of-volatility, or the edge effects arising from aggregating local staleness estimates.


We estimate spot systematic and idiosyncratic volatility in efficient price processes using local factor analysis and derive corresponding integrated volatility measures by aggregating non-overlapping local volatility proxies. We find that the volatility estimates remain unbiased, whereas estimated co-volatilities are biased due to price staleness.  By locally correcting for this bias using inverse staleness weighting, we obtain a consistent and unbiased estimator. The convergence rates of the integrated estimators are significantly faster than those of the spot estimates. This difference validates the efficiency of the aggregation process following nonlinear factor analysis. Our empirical study demonstrates that the LPCA estimator of the volatility matrix without data staleness correction results in higher out-of-sample risk in constrained portfolio allocation compared to the corrected estimator.

The remainder of this paper is organized as follows. Section \ref{sec:Price Staleness Factor Analysis} introduces the SFM, detailing the model estimation procedure and presenting the key theoretical results. Section \ref{sec:Efficient Price Volatility Estimation} describes the estimation method for efficient price volatility matrices and derives the associated theoretical properties. Section \ref{sec:Simulation} presents a simulation study that assesses the finite-sample performance of the proposed estimators. Section \ref{sec:Empirical Application} provides an empirical analysis, demonstrating the practical application of the model. Finally, Section \ref{sec:Conclusion} concludes the paper. All proofs and extra results are provided in the supplementary materials.

NOTATION. We denote by $\|\cdot\|$ and $\|\cdot\|_F$ the spectral (Euclidean) norm and the Frobenius norm, respectively. For a $d\times d$ matrix $A$ and a positive definite matrix $\Sigma$, define the weighted quadratic norm $\|A\|_{\Sigma} := d^{-1/2}\|\Sigma^{-1/2}A\Sigma^{-1/2}\|_F$. We use $a\wedge b$ and $a\vee b$ to denote $\min\{a,b\}$ and $\max\{a,b\}$. Let $\mathbf{1}_d$ be the $d$-vector of ones, and $\mathds{1}_{\{\cdot \}}$ the indicator function. Let $\lambda_1(A)\geq\lambda_2(A)\geq...\geq \lambda_{\min}(A)$ denote the ordered eigenvalues of a matrix $A$. The operator $\circ$ denotes the Hadamard product, $I_r$ is the $r\times r$ identity matrix, and $C$ denotes a generic positive constant whose value may change from line to line. For convergence, $\stackrel{P}{\longrightarrow}$ denotes convergence in probability, while $\mathcal{L}|\mathcal{F}$ and $\mathcal{L}_s|\mathcal{F}$ denote conditional weak convergence and conditional stable convergence in law given a sub-$\sigma$-field $\mathcal{F}$, respectively. For a function $f$, $f^{(i)}$ is its $i$th derivative, and $\mathcal{C}^q$ is the space of $q$-times continuously differentiable functions. We work on a probability space equipped with three information flows: (i) $(\mathcal{F}_t^{(p)})_{t\geq0}$, the natural filtration of the staleness probability process; (ii) $\mathcal{F}^{(b)}_{t_j,n}$, the $\sigma$-algebra generated by the random variables $\{b_{t_0,n},b_{t_1,n},\cdots,b_{t_j,n}\}$, forming a discrete filtration in the interval $[0,T]$; and (iii) $(\mathcal{F}_t)_{t\geq0}$, the natural filtration of the efficient price process. We further define $\mathcal{F}_{\infty}=\vee_{t>0}\mathcal{F}_t$.

\section{Price Staleness Factor Analysis}\label{sec:Price Staleness Factor Analysis}

\subsection{Price Staleness Factor Model}\label{subsec:Price Staleness}
We observe a large $d$-dimensional panel of asset log-prices, $\widetilde{Y}_{t_j}=(\widetilde{Y}_{1t_j},...,\widetilde{Y}_{dt_j})'$ sampled at equally spaced times $t_j=j\Delta_n$ for $j=0,1...,n$ over $[0,T]$, where $\Delta_n$ is the mesh and $n=\lfloor T/\Delta_n\rfloor$.\footnote{The real-world high-frequency trading data is typically asynchronous, which can be addressed using the ``\emph{previous-tick}'' synchronization scheme.} Each observed price $\widetilde{Y}_{it_j}$ either updates to the latent efficient price $Y_{it_j}$ or remains at its previous value $\widetilde{Y}_{t_{j-1}}$, depending on a Bernoulli indicator. The efficient log-price process $Y_t$ is assumed to be a $d$-dimensional It$\hat{\text{o}}$ semimartingale defined on a filtered probability space $(\Omega,\mathcal{F},(\mathcal{F}_t),\mathbb{P})$. Extending model \eqref{unimodel} to the multivariate setting gives
\vspace{-5pt}
\begin{align}\label{eq:observable price processes}
	 \widetilde{Y}_{t_j}=Y_{t_{j}}\circ(\mathbf{1}_d-B_{t_j})+\widetilde{Y}_{t_{j-1}}\circ B_{t_j},
\end{align}
where $B_{t_{j}}=\left(B_{1t_{j}},\dots,B_{dt_{j}}\right)'$ is a vector of Bernoulli random variables.

Most previous studies in the high-frequency data analysis literature have ignored the existence of price staleness (i.e., $B_{t_j}=0$ is typically assumed); c.f., \cite{mykland2009inference}, \cite{ait2017using}, \cite{kong2018systematic}, \cite{pelger2019large}, and \cite{li2024estimating}. We rewrite the Bernoulli random variable $B_{it}$ as $B_{it}=\mathds{1}_{\{b_{it}\leq p_{it}\}}$, where $\{b_{it}\}_{t\in[0,T]}$ is a collection of uniformly distributed random variables. Given the information set $\mathcal{F}^{(p)}$, the Bernoulli random variables $B_{it}$ and $B_{ms}$ are independent $\forall ~t\neq s$ or $i\neq m$. In addition, $p_t=(p_{1t}, ..., p_{dt})^{\prime}$ is modeled as a continuous-time stochastic process to capture how likely the zeros occur, which is conditionally independent of the efficient price and its volatility, given the filtration $\mathcal{F}^{(p)}$. Inspired by the generalized linear model, we define $p_{it}=\Psi(z_{it})$, where $\Psi$: $\mathbb{R}\to (0,1)$ is an increasing function in $\mathcal{C}^3$. The function $\Psi(\cdot)$ is a pre-specified link function (such as logit or probit). Our theoretical framework applies to a large class of link functions that satisfy the regularity conditions outlined in Assumption \ref{assump:staleness factor model}.3.

The process $z_{it}$ is modeled as an It$\hat{\text{o}}$ semimartingale:
\vspace{-7pt}
\begin{align*}
	z_{it} = a_i'x_{it}+\gamma_{i}'g_t, \ \ i=1,...,d,
\end{align*}
where $x_{it}$ is an $r_x$-dimensional covariate process, $a_i$ is the coefficient vector, $g_t$ is an $r_g$-dimensional continuous-time factor process independent of $\{x_{it}\}$, and $\gamma_i$ is a vector of factor loadings describing the exposure to the systematic factors. 

This two-component structure for the index $z_{it}$ is a key feature of high-dimensional staleness. The observed covariates $x_{it}$ are used to capture the effects of known, asset-specific drivers of price staleness. As demonstrated by \cite{bandi2020zeros}, variables such as trading volume are significant determinants of the probability of price updates. The latent factors $g_t$ are introduced to account for unobserved, systematic components that drive staleness co-movement across a broad set of assets. These common factors would also prevent omitted-variable bias.

We assume the processes $x_{it}$ and $g_t$ are locally bounded It$\hat{\text{o}}$ semimartingales,
\begin{align*}	
	x_{it}=x_{i0}+\int_0^t\mu_{is}^xds+\int_{0}^t\sigma_{is}^x\circ dW_{is}^x,\quad g_t=g_0+\int_0^t\mu_s^gds+\int_{0}^t\sigma_s^g\circ dW_s^g,
\end{align*}
where $W_{it}^x$ and $W_t^g$ are $r_x$-dimensional and $r_g$-dimensional Brownian motions, respectively. The coefficients $\mu_{it}^x$ and $\mu_t^g$ are progressively measurable, and $\sigma_{it}^x$ and $\sigma_t^g$ are adapted c$\grave{\text{a}}$dl$\grave{\text{a}}$g processes. There exists a uniform bound for all processes $(\mu_{is}^x,\sigma_{is}^x)$ and $(\mu_s^g,\sigma_s^g)$. Notably, we only observe the stochastic process $x_{it}$ and the Bernoulli random variables $B_{it}$, but not $p_{it}$ or $z_{it}$. This poses a challenge that the GPCA in \cite{ait2017using} and \cite{pelger2019large} and the LPCA in \cite{kong2017number,kong2018systematic}, \cite{ait2019principal}, \cite{chen2020five}, \cite{kong2023discrepancy}, and \cite{li2024estimating} are not applicable any more, because the differential form of $z_{it}$ (or $p_{it}$) is no longer observable at discrete time instances. A new method that can handle the nonstationary integral form of $z_{it}$ with continuous-time factor structure has to be invented. While it would be interesting to consider jumps in these processes, this paper does not include them in $x_{it}$ and $g_t$ due to the added complexity they introduce in our proposed MLE.\footnote{In our binary observables, the usual techniques, e.g., the truncation method in \cite{mancini2009non}, for dealing with jumps are no longer applicable.} The consideration of jumps is left for future work.

Before giving the maximum likelihood estimation method for a latent nonlinear nonstationary factor model, we give some regularity assumptions on the staleness factor model.

\begin{assumption}\label{assump:staleness factor model}
	\begin{enumerate}[label=\arabic*., leftmargin=*,rightmargin=0pt]
		\item Assume that $\|d^{-1}\Gamma'\Gamma-I_{r_g}\|\to 0$ as $d\to0$, where $\Gamma=(\gamma_1,...,\gamma_d)'$ and each $\gamma_i$ satisfies $\max_{1\leq i\leq d}\|\gamma_i\|_{F}\leq C$. Moreover, we also assume that $\sup_{t\in[0,T]}\|x_{it}\|_{F}\leq C$ and $\sup_{t\in[0,T]}\|g_{t}\|_{F}\leq C$.
		\item There exists a constant $\overline{p}\in(0,1)$ such that $\sup_{t\in[0,T]}\max_{1\leq i\leq d}p_{it}\leq \overline{p}$. Moreover, $\inf_{t\in[0,T]}\min_{1\leq i\leq d}p_{it}> 0$. 
		\item There exists a constant $C$ such that for all $z$ in $\Xi_z=\left\{z:0<\Psi(z)\leq\overline{p}\right\}$, and for $j=0,1,2$,  $|\psi^{(j)}(z)|<C$.
	\end{enumerate}
\end{assumption}

Assumption \ref{assump:staleness factor model}.1 is a strong factor condition and requires the factors to be uniformly bounded on $[0,T]$, which is standard in high-frequency factor analysis; see \cite{ait2017using}, \cite{kong2017number, kong2018systematic}, and \cite{li2024estimating}.\footnote{We assumed the strong factor condition because we established the second-order properties (the central limit theorems) of the estimated price staleness probabilities and their functionals. The strong factor condition can be relaxed to the weak one if only the first-order consistency is in need with modification of the consistency rates and more strengthened conditions such as the relation between the dimensionality and the sample size.} Assumption \ref{assump:staleness factor model}.2 requires that price staleness exists with positive probability but cannot approach one, which is mild and appears in \cite{bandi2023discontinuous}. Assumption \ref{assump:staleness factor model}.3 is a regularity condition for the link function which is satisfied by the logit and probit and many other link functions.

\begin{remark}
	Our SFM is distinguished from several classes of classic factor models, creating unique technical challenges. Unlike linear factor models for low-frequency data (e.g., \citealt{bai2003inferential,fan2013large}), our framework is inherently nonlinear and tailored for high-frequency binary observations. Compared to existing nonlinear models that operate on low-frequency data (e.g., \citealt{chen2021quantile}), our high-frequency setting requires handling nonstationary processes. Most notably, in contrast to high-frequency linear factor models based on continuous returns (e.g.,  \citealt{ait2017using,pelger2019large}), we model the probability of a discrete, latent event. This key structural difference renders standard PCA-based methods, which rely on differenced price data, inapplicable. Consequently, we must employ a novel MLE framework and develop the asymptotic theory to accommodate the model's nonlinear, nonstationary, and latent variable structure.
\end{remark}
\subsection{Estimation of the Staleness Factor Model}\label{subsec:Estimation of staleness factor}

To estimate the SFM, we employ the MLE. Define the increments of the observed covariate $x_i$ and latent factor $g$ by
\vspace{-10pt}
\begin{align*}
	\Delta x_{it_j}:=x_{it_j}-x_{it_{j-1}} \ \ \text{and} \ \
	\Delta g_{t_j}:=g_{t_j}-g_{t_{j-1}},
\end{align*}	
for $j=1,...,n$. We use the convention that $\Delta x_{it_0}:=x_{it_0}$ and $\Delta g_{t_0}:=g_{t_0}$. We next rewrite $z_{it_j}$ in the cumulative-sum representation $z_{it_j} = a_i'\sum_{l=0}^j\Delta x_{it_l}+\gamma_{i}'\sum_{l=0}^j\Delta g_{t_l}$. Since the latent index $z_{it_j}$ is unobserved and we only observe the binary outcome $B_{it_j}$, we cannot apply a PCA-type method (e.g., \citealt{ait2019principal}). Instead, we estimate the discretized factors (or their increments) jointly with the loadings by maximizing the Bernoulli likelihood. Let
\vspace{-10pt}
\begin{align*}
	 A=(a_1,...,a_d)',~~~\Gamma=(\gamma_1,...,\gamma_d)',~~~G=(g_{t_0},g_{t_1},...,g_{t_n})',~~~\Delta G=(\Delta g_{t_0},...,\Delta g_{t_n})',
\end{align*}
and $\theta_i=(a_i',\gamma_i')'$, $\Theta=(A,\Gamma)$, $u_{it}=(x_{it}',g_t')'$. The relationship between $G$ and $\Delta G$ is
$G=\varrho\Delta G$,
where $\varrho=\left(\mathds{1}_{\{j\leq i\}}\right)_{i=1,...,n+1}^{j=1,...,n+1}$ is a $(n+1)\times(n+1)$-dimensional matrix with the lower triangular and diagonal entries being 1 and others 0.

A well-known feature of factor models is that ${\gamma_i}$ and $\Delta g_{t_j}$ (or $g_{t_j}$) are not separately identifiable without an appropriate normalization. We adopt the following normalization for the SFM:
\vspace{-5pt}
\begin{align}\label{eq:identification condition 1}
	\begin{aligned}
	 \Gamma\in\mathscr{G}:=&\left\{\Gamma\in\mathbb R^{d\times r_g}:\Gamma'\Gamma/d=I_{r_{g}}\right\},\\
	 \Delta G\in\mathcal{G}:=&\left\{\Delta G\in\mathbb R^{(n+1)\times r_g}:\Delta G'\Delta G~\text{is diagonal with distinct values}\right\}.
	 \end{aligned}
\end{align}
Now, the $\mathcal{F}^{(p)}$-conditional log-likelihood is
\vspace{-5pt}
\begin{align*}
	\mathbb{L}_{d,n}(A,\Gamma,\Delta G):=\sum_{i=1}^{d}\sum_{j=0}^{n}\left\{\left(1-B_{it_j}\right)\log \left[1-\Psi(z_{it_j})\right]+B_{it_j}\log\Psi(z_{it_j})\right\},
\end{align*}
where $z_{it_j}=\left(a_i'x_{it_j}+\gamma_{i}'\sum_{l=0}^j\Delta g_{t_l}\right)$. Then the MLE $\{\hat{A},\hat{\Gamma},\Delta\hat{G}\}$ is obtained as\footnote{Ignoring the normalization in \eqref{eq:identification condition 1}, the optimization problems in $G$ and in $\Delta G$ are equivalent via invertible linear map $G=\varrho\Delta G$, hence $\hat{G}=\varrho\Delta \hat{G}$.}
\begin{align}\label{eq:log-likelihood estimator}
	(\hat{A},\hat{\Gamma},\Delta \hat{G})=\argmax_{A\in\mathbb{R}^{d\times r_{x}},\Gamma\in\mathscr{G},\Delta G\in\mathcal{G}}\mathbb{L}_{d,n}(A,\Gamma,\Delta G).
\end{align}
Unlike high-frequency PCA (global or local), our estimator does not admit a closed-form solution. This complicates both the derivation of its large-sample properties and the computation. However, as demonstrated by Theorem \ref{thm:Statistical Properties of p}, the MLE achieves the same convergence rate as the high-frequency PCA estimation (e.g., \citealt{kong2018systematic}). Let
\vspace{-5pt}
\begin{align*}
	l_{i,j}(z_{it_j})=\left\{\left(1-B_{it_j}\right)\log \left[1-\Psi(z_{it_j})\right]+B_{it_j}\log\Psi(z_{it_j})\right\},
\end{align*} 
and define $\mathbb{L}_{i,n}(\theta_i, \Delta G)=\sum_{j=0}^{n}l_{i,j}(z_{it_j})$ and $\mathbb{L}_{d,j}(\Theta, \Delta g_{t_j})=\sum_{i=1}^{d}\sum_{l=j}^nl_{i,l}(z_{it_l})$. Next, we give the computational steps.
\begin{enumerate}[label= Step \arabic*., leftmargin=*,rightmargin=0pt]
	\item Choose initial values for $\Delta G^{(0)}$ and $\Theta^{(0)}$.
	\item For each  $i=1,...,d$, given $\Delta G^{(l-1)}$, solve $\theta_i^{(l)}=\arg\max_{\theta}\mathbb{L}_{i,n}(\theta, \Delta G^{(l-1)})$. For each $j=0,1,...,n$, given $\Theta^{(l)}$, solve $\Delta g_{t_j}^{(l)}=\arg\max_{\Delta g}\mathbb{L}_{d,j}(\Theta^{(l)}, \Delta g)$.
	\item Repeat Step 2 until the criterion: $\mathbb{L}_{d,n}(\Theta^{(l^*)},\Delta G^{(l^*)})\approx\mathbb{L}_{d,n}(\Theta^{(l^*-1)},\Delta G^{(l^*-1)})$ is met for some iteration $l^*$.
	\item Normalize $\Gamma^{(l^*)}$ and $\Delta G^{(l^*)}$ to satisfy the normalization condition given in \eqref{eq:identification condition 1}. Finally, set $G^{(l^*)}=\varrho\Delta G^{(l^*)}$.
\end{enumerate}

To obtain an initial estimate, we use a local block approach to roughly estimate the staleness probability $p_{it_j}$ by $\tilde{p}_{it_j}=\bar{k}_n^{-1}\sum_{l=0}^{\bar{k}_n-1}B_{it_{j+l}}$, where $\bar{k}_n$ is a sequence of integers that satisfy $\bar{k}_n\to\infty$ and $\bar{k}_n\Delta_n\to0$. We then apply the inverse map to obtain $\tilde{z}_{it_j}=\Psi^{-1}(\tilde{p}_{it_j})$ and regress $\tilde{z}_{it_j}$ against $x_{it_j}$ for $j=0,...,n$ to get the estimate $\widetilde{a}_i$. Next, we compute the residual $\tilde{z}_{it_j}-\widetilde{a}_i'x_{it_j}$, for which we use the high-frequency PCA based on \cite{pelger2019large} to estimate $\Gamma$ and $\Delta G$. In Step 3, we set the tolerance condition as
\vspace{-5pt}
\begin{align*}
	 \frac{1}{d}\sum_{i=1}^{d}\|a_i^{(l^*)}-a_i^{(l^*-1)}\|_{F}^2+\frac{1}{nd}\|G^{(l^*)}\Gamma^{(l^*)}-G^{(l^*-1)}\Gamma^{(l^*-1)}\|_{F}^2<\varepsilon^*,
\end{align*}
for sufficiently small $\varepsilon^*>0$, e.g., $10^{-3}$. 

It is clear that our estimation steps differ from the PCA-based analytical method of \cite{ait2017using} and \cite{pelger2019large}. We employ an iterative optimization approach to obtain the estimates. In each iteration step, our algorithm is a simple optimization problem, which is consistent with the logic in \cite{chen2021quantile}. For example, if the link function is logit, each step of the iterative optimization is equivalent to a logit regression. The block-wise convexity guarantees a global optimum for each subproblem.

To mitigate the risk of not converging due to the joint non-convexity, we propose a data-driven initialization strategy based on local averaging, which provides an effective and analytical starting point. Our numerical experiments (in the Supplementary Material) confirm that this initialization performs well, that the final MLE estimates offer further improvement, and that the entire procedure robustly converges in a small number of iterations (typically fewer than 50).

To determine the number of factors consistently, we adopt \cite{pelger2019large}'s perturbed-eigenvalue ratio method, which examines the ratio of adjacent eigenvalues. We first compute the eigenvalues of $(\hat{\Gamma}\Delta \hat{G}')(\hat{\Gamma}\Delta \hat{G}')'$ and order them as $\lambda^*_1\geq\cdots\geq\lambda^*_{r_g^{\max}}$, where $r_g^{\text{max}}$ is a user-specified upper bound. After that we define perturbed eigenvalues $\hat{\lambda}^*_k=\lambda^*_k+\xi_{nd}$ where $\xi_{nd}$ is any slowly diverging sequence such that $\xi_{nd}/d\to 0$ and $\xi_{nd}\min(\sqrt{n},\sqrt{d})/d\to\infty$ (see supplementary material for details). Letting $ER_k=\hat{\lambda}^*_k/\hat{\lambda}^*_{k+1}$, we select
\vspace{-5pt}
\begin{align*}
	\hat{r}_g(\chi)=\max\{k\leq r_g^{\max}-1:ER_k>1+\chi\}, \ \mbox{for some} \ \chi>0.
\end{align*}

\subsection{Results for Staleness Factor Analysis}\label{subsec:Staleness factor model}

Let $\omega_{nd}=\min(\sqrt{n},\sqrt{d})$ and we use the infill asymptotic regime $\Delta_n\to 0$ (with $T$ fixed and $n\to \infty$) as typical in the high-frequency data analysis. Recall from Section \ref{subsec:Price Staleness} that $u_{it}=(x_{it}',g_t')$. We now introduce some more notations that pertain to the asymptotic variances. Let
$$
\Omega_u=\text{diag}\{\Omega_{u,1},...,\Omega_{u,d}\}, \ \Omega_{\gamma}=\text{diag}\{\Omega_{\gamma,1},...,\Omega_{\gamma,n+1}\}, \ \Omega_{u\gamma}=\{\Omega_{u\gamma,ij}\}_{d(r_x+r_g)\times(n+1)r_g},
$$
\vspace{-30pt}
\begin{align*}
	 &\Omega_{u,i}=\frac{1}{T}\int_0^T\frac{\psi^2(z_{it})}{\Psi(z_{it})(1-\Psi(z_{it}))}u_{it}u_{it}'dt,~~~\Omega_{\gamma,j}=\text{plim}_{d\to\infty}\frac{1}{d}\sum_{i=1}^d\frac{\psi^2(z_{it_j})}{\Psi(z_{it_j})(1-\Psi(z_{it_j}))}\gamma_i\gamma_i',\\
	 &\Omega_{u\gamma,ij}=\frac{\psi^2(z_{it_j})}{\Psi(z_{it_j})(1-\Psi(z_{it_j}))}u_{it_j}\gamma_i'.
\end{align*}

We make some assumptions about these asymptotic variances and replace $t_j$ in $\Omega_{\gamma,j}$ with $t$, renaming it $\Omega_{\gamma,t}$.
\begin{assumption}\label{assump:staleness factor model2}
	\begin{enumerate}[label=\arabic*., leftmargin=*,rightmargin=0pt]
		\item  $\sup_{t\in[0,T]}\|\frac{1}{d}\sum_{i=1}^d\frac{\psi^2(z_{it})}{(1-\Psi(z_{it}))\Psi(z_{it})}\gamma_{i}\gamma_{i}'-\Omega_{\gamma, t}\|_F=o_P(1)$ as $d\to\infty$.
		\item $\Omega_{u,i}$ and $\Omega_{\gamma,j}$ are positive definite for $1\leq i\leq d$ and $0\leq j\leq n$. $\lambda_{\max}(\Omega_u)$, $\lambda_{\max}(\Omega_{\gamma})$, $\lambda_{\max}(\Omega_u^{-1})$, $\lambda_{\max}(\Omega_{\gamma}^{-1})$, $\lambda_{\max}\left(\frac{1}{nd}\Omega_{u\gamma}'\Omega_{u}^{-1}\Omega_{u\gamma}\right)$, and $\lambda_{\max}\left(\frac{1}{nd}\Omega_{u\gamma}\Omega_{\gamma}^{-1}\Omega_{u\gamma}'\right)$ are all finite.
	\end{enumerate}
\end{assumption}

Assumption \ref{assump:staleness factor model2} provides high-level regularity and identification conditions, which are underpinned by more primitive assumptions. The uniform convergence in Assumption \ref{assump:staleness factor model2}.1 is ensured by standard moment conditions on the factor loadings $\gamma_i$ (e.g., i.i.d. draws with finite fourth-order moments) combined with the smoothness of the link function and the continuity of the latent index paths. The positive definiteness conditions in Assumption \ref{assump:staleness factor model2}.2 serve as crucial identification requirements. For $\Omega_{u,i}$, it is a standard full-rank condition precluding multicollinearity between the covariates and factors in the time series. For $\Omega_{\gamma,j}$, it is satisfied if the factor loadings exhibit sufficient cross-sectional heterogeneity, such that $\Gamma'\Gamma/d$ converges to a positive definite matrix---a direct consequence of our strong factor assumption in Assumption \ref{assump:staleness factor model}. Finally, the boundedness of the eigenvalues of these matrices and their inverses is a regularity condition implied by the bounded moments of the underlying processes (Assumption \ref{assump:staleness factor model}) and the aforementioned identification conditions, which prevent the matrices from becoming degenerate.

\begin{proposition}\label{prop:convergence rate of factor G}
 If Assumptions \ref{assump:staleness factor model} and \ref{assump:staleness factor model2} hold, and if there exists a constant $\delta^{\dag}>0$ such that $\frac{d}{n^{1+\delta^{\dag}}}=o(1)$.
	\begin{enumerate}[label=(\roman*), leftmargin=*,rightmargin=0pt]
		\item $
		 \frac{1}{\sqrt{d}}\|\hat{\Theta}-\Theta\|_F=O_{P}(\omega_{nd}^{-1})$, $\|\hat{g}_{t_j}-g_{t_j}\|=O_{P}(\omega_{nd}^{-1})$, $|\hat{\gamma}_i'\hat{g}_{t_j}-\gamma_i'g_{t_j}|=O_{P}(\omega_{nd}^{-1}).$
		\item As $\omega_{nd}\to\infty$, we have $\sum_{j=1}^n(\hat{a}_i'\Delta x_{it_j})(\hat{a}_m'\Delta x_{mt_j})=a_i'[x_i,x_m]_Ta_m+O_{P}(n^{-1/2})$, and if $n/d\to0$,
		\begin{align*}
			&\sum_{j=1}^n\Delta \hat{g}_{t_j}\Delta \hat{g}_{t_j}'=[g,g]_T+o_{P}(1),~~\sum_{j=1}^n(\hat{\gamma}_i'\Delta \hat{g}_{t_j})(\hat{\gamma}_m'\Delta \hat{g}_{t_j})=\gamma_i'[g,g]_T\gamma_m+o_{P}(1).
		\end{align*}
	\end{enumerate}
\end{proposition}

Proposition \ref{prop:convergence rate of factor G} establishes the convergence rates for the estimators and their quadratic variations. In high-frequency binary estimation, the stringent requirements on the sample size $n$ distinguish it from long-span models. Specifically, the condition $\frac{d}{n^{1+\delta^{\dag}}}=o(1)$ governs the cross-sectional maximum error for the discrete approximation of second-order moments $\Omega_u$. Estimating the quadratic variations of observable covariates is relatively straightforward. However, additional consistency conditions are required for latent factors due to the complexity of their estimation.

We now demonstrate that the estimators for the factor loadings and factors converge stably in law to mixed Gaussian distributions.\footnote{The classical results on stable convergence proposed by \cite{hall2014martingale} do not hold under the filtration $\mathcal{F}^{(b)}_{t_n,n}$, as the condition of nested filtrations is no longer satisfied. Nonetheless, this issue can be addressed using Theorem 1 and Corollary 3 from \cite{kolokolov2020statistical}.}
\begin{proposition}\label{prop:asymptotic distribution of factor G}
	Under the conditions in Proposition \ref{prop:convergence rate of factor G}, as $\omega_{nd}\to\infty$, the following pointwise convergence results hold:
	\begin{enumerate}[label=(\roman*), leftmargin=*,rightmargin=0pt]
		\item For each fixed asset $i=1,..,d$, if $\sqrt{n}/d\to0$,
		\begin{align*}
			 n^{1/2}\left(\hat{\theta}_{i}-\theta_i\right)\stackrel{\mathcal{L}_s|\mathcal{F}^{(p)}}{\longrightarrow}\mathcal{N}(0,\Omega_{u,i}^{-1}).
		\end{align*}
		\item For each fixed time point $t_j$, if $\sqrt{d}/n\to0$,
		\begin{align*}
			 d^{1/2}\left(\hat{g}_{t_j}-g_{t_j}\right)\stackrel{\mathcal{L}|\mathcal{F}^{(p)}}{\longrightarrow}\mathcal{N}(0,\Omega_{\gamma,j}^{-1}).
		\end{align*}
	\end{enumerate}
\end{proposition}

The limiting distribution of $\hat{\theta}_i$ is driven by the serial partial sums of the weighted Bernoulli variates, whereas the limiting distribution of $\hat{g}_{t_j}$ arises from their cross-sectional partial sums. Our framework inherently addresses cross-sectional dependence via the common factor structure, requiring no extra dependence assumptions for the stable convergence mode. The latent factors $g_t$ induce co-movement in the staleness probabilities $p_{it}$ across assets. Crucially, we assume that the staleness indicators $B_{it}$ are conditionally independent across $i$ given the factors and covariates. This assumption implies that all systematic cross-sectional dependence is channeled through the common factors.

Based on Propositions \ref{prop:convergence rate of factor G} and \ref{prop:asymptotic distribution of factor G}, we establish the consistency and asymptotic normality for the estimated $p_{it_j}$.

\begin{theorem}\label{thm:Statistical Properties of p}
If Assumptions \ref{assump:staleness factor model} and \ref{assump:staleness factor model2} hold, and if there exists a constant $\delta^{\dag}>0$ such that $\frac{d}{n^{1+\delta^{\dag}}}=o(1)$,
	\begin{enumerate}[label=(\roman*), leftmargin=*,rightmargin=0pt]
		\item  $\hat{p}_{it_j}-p_{it_j}=O_{P}(\omega^{-1}_{nd})$ for $i=1,...,d$.
		\item  $\omega_{nd}(\hat{p}_{it_j}-p_{it_j})/\Omega_{it_j}^{(p)}\stackrel{\mathcal{L}|\mathcal{F}^{(p)}}{\longrightarrow}\mathcal{N}_{1}$, where $\mathcal{N}_{1}$ is defined on an extension of the probability space and follows $\mathcal{N}(0,1)$ conditional on $\mathcal{F}^{(p)}$. The asymptotic variance is given by
		\begin{align}\label{eq:asympotic variance of p}
			 \Omega_{it_j}^{(p)}=\psi^2(z_{it_j})\left(\frac{\omega_{nd}^2}{n}u_{it_j}'\Omega_{u,i}^{-1}u_{it_j}+\frac{\omega_{nd}^2}{d}\gamma_i'\Omega_{\gamma,j}^{-1}\gamma_i\right).
		\end{align}
	\end{enumerate}
\end{theorem}

Theorem \ref{thm:Statistical Properties of p}(ii) manifests two notable special cases: (i) if $d/n\to0$, $\sqrt{d}(\hat{p}_{it_j}-p_{it_j})\stackrel{\mathcal{L}|\mathcal{F}^{(p)}}{\longrightarrow}\mathcal{N}\left(0,\psi^2(z_{it_j})\gamma_i'\Omega_{\gamma,j}^{-1}\gamma_i\right)$; (ii) if $n/d\to0$,  $\sqrt{n}(\hat{p}_{it_j}-p_{it_j})\stackrel{\mathcal{L}|\mathcal{F}^{(p)}}{\longrightarrow}\mathcal{N}\left(0,\psi^2(z_{it_j})u_{it_j}'\Omega_{u,i}^{-1}u_{it_j}\right)$. This is because $\hat{p}_{it_j}$ relies on the $i$th serial partial sums and $j$th cross-sectional partial sums of the Bernoulli variates.

To make the CLT feasible, one needs a consistent estimator, $\hat{\Omega}_{it_j}^{(p)}$, of the conditional variance $\Omega_{it_j}^{(p)}$ in \eqref{eq:asympotic variance of p}. Based on Proposition \ref{prop:convergence rate of factor G} and Theorem \ref{thm:Statistical Properties of p}(i), a consistent estimator of $\Omega_{it_j}^{(p)}$ can be constructed as follows (denoted by $\hat{\Omega}_{it_j}^{(p)}$):
\begin{align*}
	 \psi^2(\hat{z}_{it_j})\omega_{nd}^2\left[\hat{u}_{it_j}'\left(\sum_{j=0}^n\frac{\psi^2(\hat{z}_{it_j})\hat{u}_{it_j}\hat{u}_{it_j}'}{\Psi(\hat{z}_{it_j})\left(1-\Psi(\hat{z}_{it_j})\right)}\right)^{-1}\hat{u}_{it_j}+\hat{\gamma}_i'\left(\sum_{i=1}^d\frac{\psi^2(\hat{z}_{it_j})\hat{\gamma}_i\hat{\gamma}_i'}{\Psi(\hat{z}_{it_j})\left(1-\Psi(\hat{z}_{it_j})\right)}\right)^{-1}\hat{\gamma}_i\right],
\end{align*}
where $\hat{u}_{it_j}=(x_{it_j}',\hat{g}_{t_j}')'$ and $\hat{z}_{it_j}=\hat{a}_i^{\prime}x_{it_j}+\hat{\gamma}_i^{\prime}\hat{g}_{t_j}$. By the mode of stable convergence and since $\Omega_{it_j}^{(p)}$ is $\mathcal{F}_{\infty}^{(p)}$ measurable, we soon have the following corollary.

\begin{corollary}\label{corol:Statistical Properties of p by estimator}
	Under the conditions in Theorem \ref{thm:Statistical Properties of p},
	\begin{align*}
		 \frac{\omega_{nd}}{\sqrt{\hat{\Omega}_{it_j}^{(p)}}}(\hat{p}_{it_j}-p_{it_j})\stackrel{\mathcal{L}|\mathcal{F}^{(p)}}{\longrightarrow}\mathcal{N}(0,1),
	\end{align*}
	where $\mathcal{N}(0,1)$ is a standard normal random variable and independent of $\mathcal{F}^{(p)}$.
\end{corollary}

Besides the pointwise convergence as shown in Theorem \ref{thm:Statistical Properties of p} and Corollary \ref{corol:Statistical Properties of p by estimator}, we next introduce a global convergence result of the estimated processes in the whole time window. The integral functional of two staleness probability processes (which can naturally be generalized to the multivariate case) is useful (see Theorem \ref{thm:rate of convergence of the corrected estimators} below). Define a function $\phi$: $\Xi_p^2\to \mathbb{R}$ to be locally bounded and in $\mathcal{C}^2$, where $\Xi_p=\{p:0<p\leq\bar{p}\}$, we are interested in the following integral functional:
\vspace{-5pt}
\begin{align*}
	U_{im}(\phi):=\int_0^T\phi(p_{it},p_{mt})dt \quad \text{for} \ i\neq m.
\end{align*}
A natural estimator is
\vspace{-15pt}
\begin{align*}
	 \hat{U}_{im}^n(\Delta_n,\phi):=\Delta_n\sum_{j=0}^n\phi(\hat{p}_{it_j},\hat{p}_{mt_j}).
\end{align*}

The following theorem provides the asymptotic properties of the estimated functionals.
\begin{theorem}\label{thm:asympotic properties of functional}
 Assume that $|\partial^{j,k} \phi(x,y)|\leq C(1+|x|^{q'-j}+|y|^{q'-k})$ for $j,k=0,1,2$ and $q'\geq2$. If Assumptions \ref{assump:staleness factor model} and \ref{assump:staleness factor model2} hold, and there exists a constant $\delta^{\dag}$ such that $\frac{d}{n^{1+\delta^{\dag}}}=o(1)$, as $\min (d, n) \to \infty$,
	\begin{enumerate}[label=(\roman*), leftmargin=*,rightmargin=0pt]
		\item $\hat{U}_{im}^n(\Delta_n,\phi)\stackrel{P}{\longrightarrow}\int_0^T\phi(p_{it},p_{mt})dt$.		 
		\item If $n/d\to0$,
		$ \Delta_n^{-1/2}\left(\hat{U}_{im}^n(\Delta_n,\phi)-U_{im}(\phi)\right)\stackrel{\mathcal{L}_s|\mathcal{F}^{(p)}_{\infty}}{\longrightarrow}\frac{1}{\sqrt{T}}\left(\int_0^T\partial_1\phi(p_{it},p_{mt})u_{it}'dt\right)\Omega_{u,i}^{-1}\mathcal{N}_2$\\
		 $+\frac{1}{\sqrt{T}}\left(\int_0^T\partial_2\phi(p_{it},p_{mt})u_{mt}'dt\right)\Omega_{u,m}^{-1}\mathcal{N}_3,$\\
		where $\mathcal{N}_2$ and $\mathcal{N}_3$ are defined on an extension of the original probability space, with $\partial_1\phi(x,y)=\frac{\partial\phi(x,y)}{\partial x}$ and $\partial_2\phi(x,y)=\frac{\partial\phi(x,y)}{\partial y}$. Conditional on $\mathcal{F}^{(p)}$, the variables $\mathcal{N}_2$ and $\mathcal{N}_3$ are independent centered Gaussian random variables with covariance matrices $\Omega_{u,i}$ and $\Omega_{u,m}$, respectively.
	\end{enumerate}
\end{theorem}

To make this CLT feasible, we provide the plug-in version of Theorem \ref{thm:asympotic properties of functional}(ii).
\begin{corollary}\label{corol:Flexible CLT for functional}
	Under the conditions in Theorem \ref{thm:asympotic properties of functional},
	\begin{align*}
		 \Delta_n^{-1/2}\frac{\left(\hat{U}_{im}^n(\Delta_n,\phi)-U_{im}(\phi)\right)}{\sqrt{\widetilde{\Omega}_{u,i}+\widetilde{\Omega}_{u,m}}}\stackrel{\mathcal{L}_s|\mathcal{F}^{(p)}_{\infty}}{\longrightarrow}&\mathcal{N}(0,1),
	\end{align*}
	where ($\widetilde{\Omega}_{u,m}$ \mbox{is similarly defined})
	\begin{align*}
		 \widetilde{\Omega}_{u,i}=\frac{\Delta_n}{\sqrt{T}}\sum_{j=0}^n\partial_1\phi(\hat{p}_{it_j},\hat{p}_{mt_j})\hat{u}_{it_j}'\left(\sum_{j=0}^n\frac{\psi^2(\hat{z}_{it_j})}{\Psi(\hat{z}_{it_j})(1-\Psi(\hat{z}_{it_j}))}\hat{u}_{it_j}\hat{u}_{it_j}'\right)^{-1}\sum_{j=0}^n\partial_1\phi(\hat{p}_{it_j},\hat{p}_{mt_j})\hat{u}_{it_j}.
	\end{align*}
\end{corollary}

\begin{remark}
	Unlike the local-block approach employed by \cite{kolokolov2020statistical}, we develop our estimators of $p_{it}$ and $p_{mt}$ through MLE. Block-based methods often suffer from edge effects and nonlinear bias terms (see \citealt{jacod2013quarticity,jacod2014efficient,li2019efficient}), which sensitively depend on the chosen window size. By MLE, we eliminate these distortions tied to parameter tuning while leveraging the asymptotic efficiency of maximum-likelihood estimators.
\end{remark}

\section{Efficient Price Volatility Estimation}\label{sec:Efficient Price Volatility Estimation}

\subsection{Efficient Price Process}\label{subsec:Efficient price process}

We assume the efficient log-price process $Y$ in \eqref{eq:observable price processes}, defined on a filtered probability space $(\Omega, \mathcal{F}, (\mathcal{F}_t)_{t\geq0},P)$, follows a continuous-time factor structure of the form:
\begin{align}\label{eq:price dynamic}
	 Y_{it}=Y_{i0}+\int_{0}^t\mu_{is}ds+\sum_{l=1}^r\int_{0}^t\sigma_{is}^ldW_s^l+\int_{0}^t\sigma_{is}^*dW_{is}^*,~~~1\leq i\leq d,
\end{align}
where $\mu_i$'s, $\sigma_i^l$'s, and $\sigma_i^*$'s are locally bounded and adapted processes; $W=(W^1,\cdots,W^r)^{\prime}$ represents an $r$-dimensional standard Brownian motion; and $W^*=(W_1^*,\cdots,W_d^*)^{\prime}$ denotes a $d$-dimensional Brownian motion with correlation matrix $\rho^*=(\rho^*_{im})_{d\times d}$, independent of $W$. We impose a sparsity condition on the correlation matrix $\rho^*$ which leads to a sparse structure of the integrated idiosyncratic volatility matrix:
\begin{align*}
	\Sigma^e=(\Sigma_{im}^e)_{d\times d}=\left(\int_0^T\sigma_{is}^*\rho^*_{im}\sigma^*_{ms}ds\right)_{d\times d}.
\end{align*}

\begin{assumption}\label{assump:sparsity}
	$\rho^*\in\mathcal{I}_q(m_d):=\{\rho^*:\max_{m}\sum_{i=1}^d|\rho_{im}^*|^q$ $\leq m_d\}$
	for some $0\leq q<1$ and $m_d$ is a function of $d$. In the case $q=1$, we assume that $m_d$ is uniformly bounded in $d$. 
\end{assumption}

When $q = 0$, Assumption \ref{assump:sparsity} indicates that each asset-specific factor is correlated with at most $m_d$ assets. This type of sparsity condition on idiosyncratic correlations is standard in high-dimensional volatility modeling; see \cite{kong2018systematic}.

In matrix form, \eqref{eq:price dynamic} can be rewritten as
\vspace{-10pt}
\begin{align*}
	dY_t=\mu_tdt+\sigma_tdW_t+\sigma_t^*dW^*_t,
\end{align*}
where $Y_t=(Y_{1t},\cdots,Y_{dt})'$, $\mu_t=(\mu_{1t},\cdots,\mu_{dt})'$, $\sigma^*_t=\text{diag}(\sigma^*_{1t},\cdots,\sigma_{dt}^*)$, and $\sigma_t=(\sigma_{it}^l)_{i=1,...,d}^{l=1,...,r}$ is a $d\times r$ systematic volatility matrix.

We begin by introducing regularity assumptions for the coefficient processes of $Y$. These assumptions are standard in the literature, as seen in works such as \cite{jacod2014efficient} for univariate models, and \cite{wang2010vast}, \cite{fan2012vast}, \cite{liu2014quasi}, \cite{kim2018adaptive}, \cite{kong2018systematic}, \cite{chen2024high1}, \cite{chen2024high2}, and \cite{kong2025data} for high-dimensional It$\hat{\text{o}}$ semimartingales.

\begin{assumption}\label{assump:price coefficient bounded}
	There exists an increasing sequence of stopping times $(\tau_m)_{m\geq1}$ with $\tau_m\uparrow\infty$ almost surely, and a sequence of bounded positive constants $(\varsigma_m)_{m\geq1}$, such that for all $i=1,...,d$ and $l=1,...,r$, the following conditions hold:
	\begin{enumerate}[label=\arabic*., leftmargin=*,rightmargin=0pt]
		\item For each $m\geq1$ and all $t<\tau_m$, $|Z_t|\leq \varsigma_m$ is satisfied for $Z\in\{\mu_i, \sigma_{i}^l,\sigma_{i}^*\}$.
		\item For $Z\in\{\sigma_{i}^l,\sigma_{i}^*\}$, the following hold: $|Z_{t+s}-Z_{t}|^2\leq \varsigma_m s^{1-\epsilon}$ almost surely for some $\epsilon>0$, and $\left|E_{\mathcal{F}_{t\wedge\tau_m}}(Z_{(t+s)\wedge\tau_m}-Z_{t\wedge\tau_m})\right|+\left|E_{\mathcal{F}_{t\wedge\tau_m}}(Z_{(t+s)\wedge\tau_m}-Z_{t\wedge\tau_m})^2\right|\leq\varsigma_m s$.
	\end{enumerate}
\end{assumption}

The last regularity condition holds for $\sigma_i^l$ and $\sigma_i^*$ if they follow a Brownian It$\hat{\text{o}}$ process with locally bounded coefficient processes---a condition that can be verified using the L$\acute{\text{e}}$vy continuity theorem.

\begin{assumption}\label{assump:eigenvalues of factor component}
	We assume that $\mathrm{rank}\left(\frac{\sigma_t\sigma_t'}{d}\right)=\mathrm{rank}\left(\left(\frac{\sigma_t\sigma_t'}{d}\right)\circ\mathcal{P}_t\right)=r$.  There exists a sequence of stopping times $\tau_m\uparrow\infty$ and a sequence of positive constants $\varsigma_m^*$ such that
	\begin{align*}
		\mathrm{inf}_{0\leq t\leq \tau_m}\lambda_{r}\left(\frac{\sigma_t\sigma_t'}{d}\right)\geq\varsigma_m^*,\quad\mathrm{inf}_{0\leq t\leq \tau_m}\lambda_{r}\left(\left(\frac{\sigma_t\sigma_t'}{d}\right)\circ\mathcal{P}_t\right)\geq\varsigma_m^*,
	\end{align*}
	where $\mathcal{P}_t=\left(1-\frac{p_{it}+p_{mt}-2p_{it}p_{mt}}{1-p_{it}p_{mt}}\mathds{1}_{\{i\neq m\}}\right)_{d\times d}$ is a symmetric matrix. Furthermore, for all $t\in[0,T]$, the matrices $\sigma_t\sigma_t'/d$ and $(\sigma_t\sigma_t')\circ\mathcal{P}_t/d$ almost surely have distinct first $r$ eigenvalues, and, when sorted in decreasing order:
	\begin{align*}
		&\mathrm{inf}_{0\leq t\leq \tau_m}\min_{1\leq l\leq r-1}\left|\lambda_{l+1}\left(\frac{\sigma_t\sigma_t'}{d}\right)-\lambda_{l}\left(\frac{\sigma_t\sigma_t'}{d}\right)\right|\geq\varsigma_m^*,\\
		&\mathrm{inf}_{0\leq t\leq \tau_m}\min_{1\leq l\leq r-1}\left|\lambda_{l+1}\left(\left(\frac{\sigma_t\sigma_t'}{d}\right)\circ\mathcal{P}_t\right)-\lambda_{l}\left(\left(\frac{\sigma_t\sigma_t'}{d}\right)\circ\mathcal{P}_t\right)\right|\geq\varsigma_m^*.
	\end{align*}
	\end{assumption}

Assumption \ref{assump:eigenvalues of factor component} ensures that the leading $r$ eigenvalues are distinct and remain non-crossing over the interval $[0,T]$, thereby excluding the possibility of duplicate eigenvalues. For statistical properties of sample covariance matrix eigenvalues, see \cite{hu2019high}. The specified eigenvalue gaps in this assumption guarantee the applicability of the $\text{SIN}(\Theta)$ theorem; see \cite{fan2013large}. Moreover, this assumption implies strong factors exist, resulting in a spiked volatility matrix structure in the diffusion system. While weak factor scenarios are interesting, they fall beyond this paper's scope and are deferred for future research. Consistent rank maintenance ensures stability of the factor space.

\subsection{Estimation of Efficient Price (Co-)Volatilities}\label{subsec:Estimation of efficient price covariance}

It is not clear whether conventional volatility estimates are biased due to price staleness. To address this issue, we first briefly review the LPCA method and the estimation of systematic and idiosyncratic volatility matrices. Under the efficient price process $Y$ (see \eqref{eq:price dynamic}), the spot systematic and idiosyncratic volatility matrices are defined, respectively, as
\vspace{-5pt}
\begin{align*}	
	V_s^c:=\sigma_s\sigma_s'\quad\text{and}\quad V_s^e:=\sigma_s^*\rho^*(\sigma_s^*)'.
\end{align*}
The integrated systematic and idiosyncratic co-volatilities are
\vspace{-10pt}
\begin{align*}
	 \Sigma^c_{ij}:=\int_0^TV_{ij}^c(s)ds\quad\text{and}\quad\Sigma^e_{ij}:=\int_0^TV_{ij}^e(s)ds,
\end{align*}
respectively, where $V_{ij}^c(s)$ and $V_{ij}^e(s)$ denote the $(i,j)$th entry of $V_s^c$ and $V_s^e$, respectively.

Let $\Delta_j^nY_i=Y_{it_j}-Y_{it_{j-1}}$ and $\delta_s=(\Delta_{\lceil\frac{s}{\Delta_n}+j\rceil}^nY_i/\sqrt{\Delta_n})_{i=1,...,d}^{j=1,...,k_n}\equiv(\delta_{ij}^s)_{d\times k_n}$, where $\lceil x\rceil$ denotes the smallest integer greater than or equal to $x$. Let $\mu_s=(\mu_{it_{\lceil \frac{s}{\Delta_n}+j\rceil}})_{i=1,...,d}^{j=1,...,k_n}$, $F_s=(\Delta_{\lceil\frac{s}{\Delta_n}+j\rceil}^nW^l/\sqrt{\Delta_n})_{l=1,...,r}^{j=1,...,k_n}\equiv(F_s(1),...,F_s(k_n))$ and  $F_s^*=(\Delta^n_{\lceil\frac{s}{\Delta_n}+j\rceil}W_i^*/\sqrt{\Delta_n})_{i=1,...,d}^{j=1,...,k_n}\equiv(F^*_s(1),...,F^*_s(k_n))$. The volatility loading matrices are defined as $\sigma_s:=(\sigma_{is}^l)_{i=1,...,d}^{l=1,...,r}$ and $\sigma_s^*:=\text{diag}\{\sigma_{1s}^*,...,\sigma_{ds}^*\}$. For the window size $k_n$, we assume the following.
\begin{assumption}\label{assump:sample size}
	The sampling scheme satisfies: $k_n/\sqrt{n}=O(1)$, $\log d=o(n^{1/2-\epsilon})$, and $n/d^{2\delta'}=o(1)$ for some $\delta'\geq1$ and any $\epsilon>0$.
\end{assumption}

Following \cite{kong2018systematic}, in a local window $(s,\lceil\frac{s}{\Delta_n}\rceil\Delta_n+k_n\Delta_n)$, PCA is performed on $\frac{\delta_s'\delta_s}{dk_n}$. Specifically, $\hat{F}_s$ is the $\sqrt{k_n}$ times the eigenvector of $\frac{\delta_s'\delta_s}{dk_n}$ (with eigenvalues sorted in decreasing order) and $\hat{\sigma}_s=\frac{\delta_s\hat{F}_s'}{k_n}$. Then the estimators of $V_{im}^{c}(s)$, $V_{im}^{e}(s)$, $\Sigma_{im}^c$, and $\Sigma_{im}^e$ are, respectively, given by
\vspace{-10pt}
\begin{align}\label{eq:volatility estimators}
	\begin{aligned}
		 &\hat{V}_{im}^c(s)=\hat{\sigma}'_{is}\hat{\sigma}_{ms},~~~~~~~~\hat{V}^e_{ii}(s)=\frac{1}{k_n}\sum_{j=1}^{k_n}(\delta_{ij}^s)^2-\hat{V}_{ii}^c(s),\\
		 &\hat{V}^e_{im}(s)=\frac{1}{k_n}\sum_{j=1}^{k_n}(\delta_{ij}^s-\hat{\sigma}_{is}'\hat{F}_s(j))(\delta_{mj}^s-\hat{\sigma}_{ms}'\hat{F}_s(j))\quad\text{for}\ i\neq m, \\
		 &\hat{\Sigma}_{im}^c=k_n\Delta_n\sum_{k=1}^{\lfloor n/k_n\rfloor}\hat{V}_{im}^c(t_{(k-1)k_n}),~~~	 \hat{\Sigma}_{im}^e=k_n\Delta_n\sum_{k=1}^{\lfloor n/k_n\rfloor}\hat{V}_{im}^e(t_{(k-1)k_n}).
	\end{aligned}
\end{align}

For this low-rank plus sparse setting, we use the Principal Orthogonal complEment Thresholding (POET) method given in \cite{fan2013large} and \cite{kong2018systematic}. Taking the spot idiosyncratic volatility ($\hat{V}_s^{e\mathcal{T}}=(\hat{V}_{im}^{e\mathcal{T}}(s))_{d\times d}$) as an example, its $(i,m)$th entry is given by $\hat{V}_{im}^{e\mathcal{T}}(s)=\hat{V}_{ii}^{e}(s)$ for the diagonal entries ($i=m$) and by $\hat{V}_{im}^{e\mathcal{T}}(s)=s_{im}(\hat{V}_{im}^{e}(s))$ for the off-diagonal entries ($i\neq m$), where $s_{im}(\cdot)$ is a generalized shrinkage function given in \cite{fan2013large}. The integrated idiosyncratic volatility is treated analogously and is denoted as $\hat{\Sigma}^{e\mathcal{T}}=(\hat{\Sigma}_{im}^{e\mathcal{T}})_{d\times d}$. In addition, $\tau_{im}$ is an entry-dependent threshold, which is $\tau_{im}=C\varphi_{nd}\sqrt{\hat{\hslash}_{im}}$ for spot volatilities and $\tau_{im}=C\widetilde{\varphi}_{nd}\sqrt{\hat{\hbar}_{im}}$ for integrated volatilities (see Theorem \ref{thm:rate of convergence of the volatility matrices} for $\widetilde{\varphi}_{nd}$ and $\varphi_{nd}$).\footnote{Note that $\hat{\hslash}_{im}$ and $\hat{\hbar}_{im}$ are chosen similarly to \cite{fan2013large}, and we choose $\hat{\hslash}_{im}=\frac{1}{k_n}\sum_{j=1}^{k_n}[(\delta_{ij}^s-\hat{\sigma}_{is}'\hat{F}_s(j))(\delta_{mj}^s-\hat{\sigma}_{ms}'\hat{F}_s(j))-\hat{V}_{im}^e(s)]^2$ and $\hat{\hbar}_{im}=k_n\Delta_n\sum_{k=1}^{\lfloor n/k_n\rfloor}[\hat{V}_{im}^e(t_{(k-1)k_n})-\hat{\Sigma}_{im}^e]^2$.} Consequently, our factor-based estimators of the total (systematic plus idiosyncratic) spot and integrated volatility matrices are, respectively,
\vspace{-5pt}
\begin{align*}
	 \hat{V}_s=\hat{V}_s^c+\hat{V}_s^{e\mathcal{T}}~~~\text{and}~~~\hat{\Sigma}=\hat{\Sigma}^c+\hat{\Sigma}^{e\mathcal{T}}.
\end{align*}

If staleness occurs, we observe $\widetilde{Y}$, and we denote $\widetilde{\delta}_s=(\Delta_{\lceil\frac{s}{\Delta_n}+j\rceil}^n\widetilde{Y}_i/\sqrt{\Delta_n})_{i=1,...,d}^{j=1,...,k_n}$. In a local window $(s,\lceil\frac{s}{\Delta_n}\rceil\Delta_n+k_n\Delta_n)$, we denote $B_{i\lceil\frac{s}{\Delta_n}\rceil+j}=B_{si}(j)=B_s(i,j)$,
\vspace{-5pt}
\begin{align*}
	 \alpha_{s,jl}^{(i)}=(1-B_s(i,j))\prod_{k=1}^lB_{s}(i,j-k)~~\text{for}~~l\geq1,~~\text{and}~~\alpha_{s,j0}^{(i)}=(1-B_s(i,j)).
\end{align*}
Here $\alpha_{s,jl}^{(i)}$ is the indicator that the most recent price update of asset $i$ within the window $(t_{j-l},t_j]$ occurs at $t_{j-l}$. Thus, we can express $\widetilde{\delta}_s$ in the following form.
\vspace{-5pt}
\begin{align*}
	 \widetilde{\delta}^s_{ij}=\Delta_{\lceil\frac{s}{\Delta_n}+j\rceil}^n\widetilde{Y}_i/\sqrt{\Delta_n}=\sum_{l=0}^{j-1}\alpha_{s,jl}^{(i)}\Delta_{\lceil\frac{s}{\Delta_n}+j-l\rceil}^nY_i/\sqrt{\Delta_n}=\sum_{l=1}^j\alpha_{s,j(j-l)}^{(i)}\Delta_{\lceil\frac{s}{\Delta_n}+l\rceil}^nY_i/\sqrt{\Delta_n},
\end{align*}
and the relationship between $\widetilde{\delta}_s$ and $\delta_s$ is $\widetilde{\delta}^s_{ij}=\sum_{l=1}^j\alpha_{s,j(j-l)}^{(i)}\delta^s_{il}$. Interestingly, introducing price staleness in our model is akin to incorporating factor lags; however, our model adds complexity by utilizing random coefficients. To determine the number of factors, $r$, we use an information-criterion approach, minimizing the aggregated mean squared residual error with a penalty, as outlined in \cite{kong2017number}. In the theoretical analysis of the next section, we assume the number of factors $r$ is known. When estimating the number of factors using \cite{kong2017number}'s method, the concentration inequality yields an $o(1)$ term.

\subsection{Results of Estimating the Efficient Price (Co-)Volatilities}\label{subsec:Recovery effective price volatility matrices}

Our first result below demonstrates that ignoring the price staleness introduces bias in estimating the co-volatilities.

\begin{theorem}\label{thm:rate of convergence of the volatility matrices}
	Suppose Assumptions \ref{assump:staleness factor model}--\ref{assump:sample size} hold, $\max_{m\leq d}\sum_{i=1}^d|\rho^*_{im}|/\sqrt{d}<C$, $\lambda_{\max}(\rho^*\circ\mathcal{P}_s)<C$ for some positive constant $C$.
	\begin{enumerate}[label=(\roman*), leftmargin=*,rightmargin=0pt]
		\item The systematic (co-)volatilities satisfy
		\begin{align*}
			 \hat{V}_{im}^c(s)-\left(1-\frac{p_{is}+p_{ms}-2p_{is}p_{ms}}{1-p_{is}p_{ms}}\mathds{1}_{\{i\neq m\}}\right)\sigma_{is}'\sigma_{ms}=&O_{P}\left(\frac{1}{d\wedge n^{1/4}}\right),\\
			 \hat{\Sigma}_{im}^c-\int_0^T\left(1-\frac{p_{is}+p_{ms}-2p_{is}p_{ms}}{1-p_{is}p_{ms}}\mathds{1}_{\{i\neq m\}}\right)\sigma_{is}'\sigma_{ms}ds=&O_{P}\left(\frac{1}{d\wedge n^{1/2}}\right).
		\end{align*}
		\item The idiosyncratic volatility matrices satisfy
		\begin{align*}
			 P\left(\sup_{\rho^*\in\mathcal{I}_q(m_d)}\|\hat{V}^{e\mathcal{T}}_{s}-V^{e,(p)}_{s})\|\leq C_qm_d\varphi_{nd}^{1-q}\right)=&1-O(d^{-\delta'}n^{1/2}+d^{-\delta'/2}+d^{1-\delta'}n^{1-\delta'/2}),\\
			 P\left(\sup_{\rho^*\in\mathcal{I}_q(m_d)}\|\hat{\Sigma}^{e\mathcal{T}}-\Sigma^{e,(p)})\|\leq C_qm_d\widetilde{\varphi}_{nd}^{1-q}\right)=&1-O(d^{-\delta'}n^{1/2}+d^{-\delta'/2}+d^{1-\delta'}n^{1-\delta'/2}),
		\end{align*}
		for some constant $C_q$, where $\varphi_{nd}=\frac{1}{\sqrt{d}}+\frac{\sqrt{\log d}}{n^{1/4}}$, $\widetilde{\varphi}_{nd}=\frac{1}{\sqrt{d}}+\frac{\sqrt{\log d}}{\sqrt{n}}$, $V^{e,(p)}_{s}=V^{e}_{s}\circ\mathcal{P}_s$, and $\Sigma^{e,(p)}=\int_0^TV^{e,(p)}_{s}ds$.
	\end{enumerate}
\end{theorem}

The process $p$ does not introduce bias in the estimates of either spot or integrated systematic volatilities ($i=m$), but it does bias the estimates of co-volatilities ($i\neq m$). Notably, our convergence rates match those for efficient price volatility estimates established in \cite{kong2018systematic}. Furthermore, we find that the $(i,m)$th entry of $\mathcal{P}_s$ equals zero if either $p_{is}$ or $p_{ms}$ attains a value of 1. In such cases, recovering the effective price co-volatility matrix is challenging, which is avoided by Assumption \ref{assump:staleness factor model}.2.

Theorem \ref{thm:rate of convergence of the volatility matrices}(ii) shows that the thresholding estimates of sparse spot and integrated idiosyncratic volatility matrices converge at rates $m_d\varphi_{nd}^{1-q}$ and $m_d\widetilde{\varphi}_{nd}^{1-q}$, respectively. Note that $V^{e,(p)}_{s}$ and $\Sigma^{e,(p)}$ are influenced by $\mathcal{P}_s$, indicating that price staleness affects both systematic and idiosyncratic co-volatilities.

\begin{remark}
	When the staleness is present, the co-volatility is underestimated. Indeed, the bias factor   $1-\frac{p_{is}+p_{ms}-2p_{is}p_{ms}}{1-p_{is}p_{ms}}$ is equal to $\frac{(1-p_{is})(1-p_{ms})}{1-p_{is}p_{ms}}$ which lies in $(0, 1)$. When either $p_{is}$ or $p_{ms}$ tends to 1, the bias factor tends toward 0; when both tend to 0, it tends to 1.
\end{remark}

In cases with highly spiked eigenvalues, covariance matrices cannot be consistently estimated in the spectral norm, but they can be accurately estimated in terms of the relative errors, as discussed by \cite{fan2013large}. Specifically, we consider the relative error matrix $V_s^{-1/2}\hat{V}_{s}V_s^{-1/2}-I_d$, measured by its normalized Frobenius norm $d^{-1/2}\|V_s^{-1/2}\hat{V}_{s}V_s^{-1/2}-I_d\|_F=:\|\hat{V}_s-V_s\|_{V_s}$. The following theorem summarizes the convergence results of the estimated total volatility matrix and its inverse.

\begin{theorem}\label{thm:rate of convergence of the price volatility matrices}
	Assume the conditions in Theorem \ref{thm:rate of convergence of the volatility matrices} hold.
	\begin{enumerate}[label=(\roman*), leftmargin=*,rightmargin=0pt]
		\item Let $\varphi_{nd}=\frac{1}{\sqrt{d}}+\frac{\sqrt{\log d}}{n^{1/4}}$, for some positive constant $C_q$,
		\begin{align*} &P\left(\sup_{\rho^*\in\mathcal{I}_q(m_d)}\|\hat{V}_{s}-V^{(p)}_{s}\|_{V^{(p)}_{s}}\leq C_q\left(m_d\varphi_{nd}^{1-q}+\frac{1}{d^{1/4}}+\frac{\sqrt{\log d}}{n^{(1-\epsilon)/4}}\right)\right)\\
			 &=1-O(d^{-\delta'}n^{1/2}+d^{-\delta'/2}+d^{1-\delta'}n^{1-\delta'/2}).
		\end{align*}
		\item If $m_d\varphi_{nd}^{1-q}=o(1)$, $d^{-\delta'}n^{1/2}+d^{1-\delta'}n^{1-\delta'/2}=o(1)$, $\inf_{s\in[0,T]}\min_{1\leq i\leq d}|\sigma_{is}^*|>c^{-1}$ and $c^{-1}\leq\lambda_{\min}(\rho^*\circ\mathcal{P}_s)\leq\lambda_{\max}(\rho^*\circ\mathcal{P}_s)\leq c$ for some positive constant $c$,
		\begin{align*} \|(\hat{V}_{s})^{-1}-(V_s^{(p)})^{-1}\|=O_{P}\left(m_d\varphi_{nd}^{1-q}+\frac{1}{\sqrt{d}}+\frac{\sqrt{\log d}}{n^{1/4}}\right).
		\end{align*}
	\end{enumerate}
\end{theorem}

In Theorem \ref{thm:rate of convergence of the price volatility matrices}, the term $d^{-1/4}+\frac{\sqrt{\log d}}{n^{(1-\epsilon)/4}}=d^{-1/4}+o(n^{-\epsilon/4})$ tends to zero based on Assumption \ref{assump:sample size}. Theorem \ref{thm:rate of convergence of the price volatility matrices} indicates that our volatility (precision) matrix estimate is not consistent with the volatility (precision) matrix of the efficient price (i.e., $V_s^{-1}$) in the presence of price staleness. A straightforward inverse probability weighting correction for $i \neq m$ is
\vspace{-15pt}
\begin{align*}
	 &\hat{V}_{im}^{c\star}(s):=\hat{V}_{im}^c(s)\phi(\hat{p}_{is},\hat{p}_{ms})^{-1},~~~\hat{V}_{im}^{e\star}(s):=\hat{V}_{im}^e(s)\phi(\hat{p}_{is},\hat{p}_{ms})^{-1},\\
	 &\hat{\Sigma}_{im}^{c\star}:=k_n\Delta_n\sum_{k=1}^{\lfloor n/k_n\rfloor}\hat{V}_{im}^c(t_{(k-1)k_n})\phi(\hat{p}_{it_{(k-1)k_n}},\hat{p}_{mt_{(k-1)k_n}})^{-1},\\
	 &\hat{\Sigma}_{im}^{e\star}:=k_n\Delta_n\sum_{k=1}^{\lfloor n/k_n\rfloor}\hat{V}_{im}^e(t_{(k-1)k_n})\phi(\hat{p}_{it_{(k-1)k_n}},\hat{p}_{mt_{(k-1)k_n}})^{-1},
\end{align*}
where $\hat{V}_{im}^c(s)$ and $\hat{V}_{im}^e(s)$ are given in \eqref{eq:volatility estimators}, $\hat{p}_{is}$ and $\hat{p}_{ms}$ are the maximum likelihood estimators in \eqref{eq:log-likelihood estimator}, and $\phi(x,y)=\frac{(1-x)(1-y)}{1-xy}$. Similarly, the idiosyncratic volatility matrix estimators can be corrected by thresholding the matrices ($\hat{V}_{im}^{e\star}(s)$ and ($\hat{\Sigma}_{im}^{e\star}$), and denoted by $\hat{V}_{s}^{e\star\mathcal{T}}$ (spot) and $\hat{\Sigma}^{e\star\mathcal{T}}$ (integrated), respectively. Define
\vspace{-5pt}
\begin{align*}
	 \hat{V}_s^{\star}=\hat{V}_s^{c{\star}}+\hat{V}_s^{e\star\mathcal{T}}~~~\text{and}~~~\hat{\Sigma}^{\star}=\hat{\Sigma}^{c{\star}}+\hat{\Sigma}^{e\star\mathcal{T}}.
\end{align*}

Like most of existing volatility matrix estimators, we cannot guarantee the positiveness of these estimators for finite sample. We use an intuitive appealing projection method in \cite{fan2012vast} that forces negative eigenvalues to be non-negative. The next theorem gives the convergence rates of the bias-corrected estimators of the systematic and idiosyncratic volatilities.
\begin{theorem}\label{thm:rate of convergence of the corrected estimators}
	Under the conditions of Theorem \ref{thm:rate of convergence of the volatility matrices} and the additional restriction $\lambda_{\max}(\rho^*)<C$ for some positive constant $C$, the following results hold:
	\begin{enumerate}[label=(\roman*), leftmargin=*,rightmargin=0pt]
		\item For systematic co-volatilities with $i\neq m$,
		\begin{align*}
			 \hat{V}_{im}^{c\star}(s)-\sigma_{is}'\sigma_{ms}=O_{P}\left(\frac{1}{d^{1/2}\wedge n^{1/4}}\right),\quad
			 \hat{\Sigma}_{im}^{c\star}-\int_0^T\sigma_{is}'\sigma_{ms}ds=O_{P}\left(\frac{1}{d^{1/2}\wedge n^{1/2}}\right).
		\end{align*}
		\item For idiosyncratic volatility matrices, assume there exist constants $\delta^{\dag}$, $\delta^{\ddag}$, and $\delta^{\S}$ such that $\frac{d}{n^{1+\delta^{\dag}}}+ \frac{n}{d^{2-\delta^{\ddag}}\log d} + \frac{d}{n^{2-\delta^{\S}}\log n} = o(1)$. Then, for some constant $C_q$,
		\begin{align*}
			 P\left(\sup_{\rho^*\in\mathcal{I}_q(m_d)}\|\hat{V}^{e\star\mathcal{T}}_{s}-V^{e}_{s}\|\leq C_qm_d\mathring{\varphi}_{nd}^{1-q}\right)=&1-O(d^{-\delta'}n^{1/2}+d^{-\delta'/2}+d^{1-\delta'}n^{1-\delta'/2}),\\
			 P\left(\sup_{\rho^*\in\mathcal{I}_q(m_d)}\|\hat{\Sigma}^{e\star\mathcal{T}}-\Sigma^{e}\|\leq C_qm_d\breve{\varphi}_{nd}^{1-q}\right)=&1-O(d^{-\delta'}n^{1/2}+d^{-\delta'/2}+d^{1-\delta'}n^{1-\delta'/2}),
		\end{align*}
		where $\mathring{\varphi}_{nd}=\frac{\sqrt{\log n}}{d^{1/2}}+\frac{\sqrt{\log d}}{n^{1/4}}$ and $\breve{\varphi}_{nd}=\frac{\sqrt{\log n}}{d^{1/2}}+\frac{\sqrt{\log d}}{n^{1/2}}$.
	\end{enumerate}
\end{theorem}

After applying the correction, the spot systematic volatility achieves a convergence rate of $d^{1/2} \wedge n^{1/4}$, while the integrated systematic volatility attains $d^{1/2} \wedge n^{1/2}$. Both estimates are asymptotically unbiased and thus robust to the data staleness, which is also true for the estimated total volatility matrix and its inverse. Note that the convergence rate in Theorem \ref{thm:rate of convergence of the volatility matrices}(i) has decreased from $d^{-1}$ to $d^{-1/2}$ in Theorem \ref{thm:rate of convergence of the corrected estimators}(i) due to the inclusion of the estimation errors of price staleness probabilities.
\begin{theorem}\label{thm:rate of convergence of the corrected price volatility matrices}
	Assume the conditions in Theorem \ref{thm:rate of convergence of the corrected estimators} hold.
	\begin{enumerate}[label=(\roman*), leftmargin=*,rightmargin=0pt]
		\item Let $\mathring{\varphi}_{nd}=\frac{\sqrt{\log n}}{d^{1/2}}+\frac{\sqrt{\log d}}{n^{1/4}}$. For some positive constant $C_q$,
		\begin{align*}
			 &P\left(\sup_{\rho^*\in\mathcal{I}_q(m_d)}\|\hat{V}_{s}^{\star}-V_{s}\|_{V_{s}}\leq C_q\left(m_d\mathring{\varphi}_{nd}^{1-q}+\frac{1}{d^{1/4}}+\frac{\sqrt{\log d}}{n^{(1-\epsilon)/4}}+\sqrt{\frac{\log n}{d}}\right)\right)\\
			 &=1-O(d^{-\delta'}n^{1/2}+d^{-\delta'/2}+d^{1-\delta'}n^{1-\delta'/2}).
		\end{align*}
		\item If $m_d\mathring{\varphi}_{nd}^{1-q}=o(1)$, $d^{-\delta'}n^{1/2}+d^{1-\delta'}n^{1-\delta'/2}=o(1)$, $\inf_{s\in[0,T]}\min_{1\leq i\leq d}|\sigma_{is}^*|>c^{-1}$ and $c^{-1}\leq\lambda_{\min}(\rho^*)\leq\lambda_{\max}(\rho^*)\leq c$ for some positive constant $c$,
		\vspace{-5pt}
		\begin{align*}
			 \|(\hat{V}_{s}^{\star})^{-1}-(V_s)^{-1}\|=O_P\left(m_d\mathring{\varphi}_{nd}^{1-q}+\frac{\sqrt{\log n}}{d^{1/2}}+\frac{1}{d^{1/4}}+\frac{\sqrt{\log d}}{n^{1/4}}\right).
		\end{align*}
	\end{enumerate}
\end{theorem}

In both bounds, the $\sqrt{\log n/d}$ term originates from the estimation of staleness probability. In other words, incorporating staleness probability estimation brings these additional $\sqrt{\log n/d}$ terms into the overall error bounds.

\section{Simulation}\label{sec:Simulation}

\subsection{Simulation Design}\label{subsec:Simulation experiment setting}

We generate one-minute (or five-minute) sampling over a 6.5-hr trading day from model (\ref{eq:observable price processes}), where the Bernoulli variates $B_{ij}$ are generated in the following steps.
\begin{enumerate}[label=Step \arabic*., leftmargin=*,rightmargin=0pt]
	\item Generate i.i.d. uniform random variables $\{b_{ij}\}_{i=1,\cdots,d;j=1,\cdots,n}\sim U(0,1)$.
	\item We generate the latent variable $z_{it_j}=a_i'x_{it_j}+\gamma_{i}'g_{t_j}$ and then set $p_{it_j}=\Psi(z_{it_j})$, where $\Psi$ is either the probit or logit link function. The loading coefficients $a_i$ and $\gamma_i$ are drawn independently across $i$. Elements of $a_i$ are from $U(0,1.5)$ and elements of $\gamma_i$ from $\mathcal{N}(0,1)$. The covariate vector $x$ and factor vector $g$ follow mean-reverting processes:
	\vspace{-10pt}
	\begin{align*}
		dx_{it}=&\kappa_{x}\circ(\mu_x-x_{it})d t+\sigma_{x} \circ d W_{i t}^{x} ,\ \
		dg_{t}=\kappa_{g}\circ(\mu_g-g_t)d t+\sigma_g \circ d W_{t}^g,
	\end{align*}
    with parameters $\kappa_x=(50,50)'$, $r_x=(5,5)'$, $\kappa_g=(10,15)'$, $\sigma_g=(1,1)'$, and mean vectors $\mu_g=(0,0)'$ and $\mu_x=p_{\text{level}}\times(0.4,0.4)'$. The long-term staleness probability is set to 0.05, yielding $p_{\text{level}}=\log(0.05/0.95)$ for the logit link. Initial values are set based on an initial probability of 0.5, with $p_{\text{initial}}=0$, such that $x_{i0}=p_{\text{initial}}\times(0.6,0.4)$ and $g_{0}=(0,0)'$.
	\item Generate Bernoulli variates from $B_{ij}=\mathds{1}_{\{b_{ij}\leq p_{it_j}\}}$.
\end{enumerate}

We model the efficient price process $Y$ following the three-factor ($r=3$) setup of \cite{kong2018systematic}, where the systematic spot volatility follows a square-root process.
\vspace{-10pt}
\begin{align*}
	d\left(\sigma_{i t}^l\right)^2=c_{l i}\left(a_{l i}-\left(\sigma_{i t}^l\right)^2\right) d t+\sigma_{l i}^0 \sigma_{i t}^l d W_{i t}^\sigma, \quad l=1, \ldots, r .
\end{align*}

We set $a_{1 i}=0.5+i / d, a_{2 i}=0.75+i / d, a_{3 i}=0.6+i / d, c_{1 i}=0.03+i / (100 d)$, $c_{2 i}=0.05+i / (100 d), c_{3 i}=0.08+i / (100 d), \sigma_{1 i}^0=0.15+i / (10 d), \sigma_{2 i}^0=\sigma_{3 i}^0=0.2+$ $i / (10 d)$. The specific volatility process follows the stochastic differential equation,
\vspace{-5pt}
\begin{align*}
	d\left(\sigma_{i t}^*\right)^2=\left(0.08+\frac{i}{100d}\right)\left(0.25+\frac{i}{d}-\left(\sigma_{i t}^*\right)^2\right) d t+\left(0.2+\frac{i}{10d}\right) \sigma_{i t}^* d W_{i t}^{\sigma *} .
\end{align*}
We set the initial values to $(\sigma_{i 0}^1,\sigma_{i 0}^2,\sigma_{i 0}^3)=(\sqrt{0.06},\sqrt{0.04},\sqrt{0.08})$ and $\sigma_{i 0}^*=\sqrt{0.03}$.

As in \cite{jacod2014efficient} and \cite{kong2018systematic}, we simulate the efficient prices from
\vspace{-30pt}
\begin{align*}
	dY_{is}=\sigma_{is}^1dW_s^1+\cdots+\sigma_{is}^rdW_s^r+\sigma_{is}^*dW_{is}^*,
\end{align*}
where $\{W_s^l\}_{l=1}^r$ and $\{W_{is}^*\}_{i=1}^{d}$ are independent Brownian motions. Moreover, for each $i$, the Brownian motions $W_{i s}^\sigma$,  $W_{i s}^{\sigma *}$, $\{W_s^l\}_{l=1}^r$, and $W_{is}^*$ are mutually independent. The correlation matrix $\rho^*$ is a block diagonal matrix with each block being $\left(\rho_{ij}^*=0.6^{|i-j|}\right)_{10\times10}$. This setting is similar to that in \cite{ait2017using} and \cite{kong2023discrepancy}. All simulations are repeated 100 times for dimensions $d\in\{50,100,150\}$, assuming the number of factors is known. We examine two scenarios over a 3-day period: one-minute data ($n=1170$) and five-minute data ($n=234$). For both, the local window size is set to $k_n\approx\sqrt{n}$ (specifically, $k_n=30$ for $n=1170$ and $k_n=15$ for $n=234$), which partitions the data into 39 and 15 blocks, respectively.

The proposed procedure is computationally efficient. For a panel dataset with $n=234$, the computation times are approximately 5, 6, and 9 seconds for assets with $d=50$, $100$, and $150$, respectively. All computations are implemented in MATLAB R2020b on a desktop computer equipped with an Intel Core i7 3.80GHz CPU.

\subsection{Simulation Results}\label{Simulation results}
To evaluate the accuracy of SFM estimates, the staleness-induced bias in volatility matrices, and the effectiveness of our correction, we report the simulation results below. To save space, additional simulation details and results are provided in the supplementary materials. Specifically, the supplementary materials report (i) a comparison of several probability paths of true staleness, MLE estimates, and local block estimates; and (ii) the finite-sample approximations to the limit distributions in Corollaries \ref{corol:Statistical Properties of p by estimator} and \ref{corol:Flexible CLT for functional}. 

\begin{table}[htbp]
	\footnotesize
	\centering
	\caption{Simulation results for staleness factor model and volatility matrix estimation.}
	\label{tab:Numerical simulation results} 
	\setlength{\tabcolsep}{3pt}
	\renewcommand{\arraystretch}{0.75}
	\begin{tabular*}{\textwidth}{@{\extracolsep{\fill}} l *{12}{c}} 
		\toprule
		& \multicolumn{6}{c}{Logit type} & \multicolumn{6}{c}{Probit type} \\
		\cmidrule(lr){2-7} \cmidrule(lr){8-13}
		& \multicolumn{3}{c}{5 minute} & \multicolumn{3}{c}{1 minute} & \multicolumn{3}{c}{5 minute} & \multicolumn{3}{c}{1 minute} \\
		\cmidrule(lr){2-4} \cmidrule(lr){5-7} \cmidrule(lr){8-10} \cmidrule(lr){11-13}
		$d$ & 50 & 100 & 150 & 50 & 100 & 150 & 50 & 100 & 150 & 50 & 100 & 150 \\
		\midrule
		\multicolumn{13}{@{}c}{\textit{Panel A: Staleness factor model results}} \\
		\midrule
		No.Factor & 1.735 & 1.862 & 1.926 & 1.837 & 1.914 & 1.998 & 2.257 & 2.165 & 2.052 & 2.174 & 2.051 & 2.000 \\
		RMSE$_z^{\text{LB}}$ & 0.610 & 0.597 & 0.589 & 0.202 & 0.201 & 0.199 & 0.379 & 0.377 & 0.375 & 0.186 & 0.184 & 0.183 \\
		RMSE$_z^{\text{MLE}}$ & 0.241 & 0.243 & 0.238 & 0.105 & 0.104 & 0.104 & 0.213 & 0.207 & 0.210 & 0.082 & 0.085 & 0.083 \\
		RMSE$_p^{\text{LB}}$ & 0.062 & 0.062 & 0.061 & 0.024 & 0.024 & 0.024 & 0.052 & 0.052 & 0.052 & 0.024 & 0.025 & 0.024 \\
		RMSE$_p^{\text{MLE}}$ & 0.040 & 0.040 & 0.040 & 0.018 & 0.018 & 0.018 & 0.035 & 0.035 & 0.035 & 0.015 & 0.016 & 0.016 \\
		\midrule
		\multicolumn{13}{@{}c}{\textit{Panel B: Without staleness (only efficient price)}} \\
		\midrule
		$\|\hat{V}_t-V_t\|$ & 1.694  & 3.721  & 5.299  & 1.196  & 2.399  & 3.716  & 1.694  & 3.721  & 5.299  & 1.196  & 2.399  & 3.716  \\
		$\|\hat{\Sigma}-\Sigma\|$ & 0.799  & 1.761  & 2.521  & 0.356  & 0.777  & 1.290  & 0.799  & 1.761  & 2.521  &  0.356  & 0.777  & 1.290  \\
		\midrule
		\multicolumn{13}{@{}c}{\textit{Panel C: With staleness \& no correction}} \\
		\midrule
		$\|\hat{V}_t-V_t^{(p)}\|$ & 2.174  & 4.394  & 6.349  & 1.941  & 3.866  & 5.808  & 2.814  & 5.705  & 8.329  & 2.682  & 5.352  & 8.075  \\
		$\|\hat{V}_t-V_t\|$ & 4.516  & 9.074  & 13.980  & 4.376  & 8.867  & 13.387  & 3.928  & 7.888  & 12.092  & 3.717  & 7.542  & 11.388  \\
		$\|\hat{\Sigma}-\Sigma^{(p)}\|$ & 0.689  & 1.424  & 2.020  & 0.316  & 0.654  & 1.017  & 0.721  & 1.534  & 2.214  & 0.340  & 0.702  & 1.122  \\
		$\|\hat{\Sigma}-\Sigma\|$ & 4.378  & 8.793  & 13.592  & 4.310  & 8.727  & 13.175  & 3.654  & 7.329  & 11.350  & 3.557  & 7.219  & 10.886  \\
		\midrule
		\multicolumn{13}{@{}c}{\textit{Panel D: With staleness \& correction}} \\
		\midrule
		$\|\hat{V}_t^{\text{O}}-V_t\|$ & 2.896  & 6.045  & 8.680  & 1.985  & 4.015  & 6.005  & 2.717  & 5.845  & 8.381  & 1.946  & 3.900  & 5.884  \\
		$\|\hat{V}_t^{\text{LB}}-V_t\|$ & 3.138  & 6.506  & 9.330  &  2.022  & 4.096  & 6.126  & 2.940  & 6.397  & 9.343  & 2.005  & 4.020  & 6.046  \\
		$\|\hat{V}_t^{\star}-V_t\|$ & 2.961  & 6.184  & 8.857  &  2.001  & 4.049  & 6.054  &  2.784  & 5.998  & 8.573  & 1.988  & 3.977  & 5.984  \\
		$\|\hat{\Sigma}^{\text{O}}-\Sigma\|$ & 1.407  & 2.906  & 4.349  & 0.652  & 1.353  & 2.101  & 1.353  & 2.889  & 4.361  & 0.650  & 1.345  & 2.108  \\
		$\|\hat{\Sigma}^{\text{LB}}-\Sigma\|$ & 1.566  & 3.274  & 4.853  &  0.735  & 1.532  & 2.361  & 1.544  & 3.365  & 5.189  & 0.770  & 1.580  & 2.448  \\
		$\|\hat{\Sigma}^{\star}-\Sigma\|$ & 1.436  & 2.994  & 4.408  & 0.665  & 1.387  & 2.131  & 1.377  & 2.966  & 4.436  & 0.667  & 1.379  & 2.140  \\
		\bottomrule
	\end{tabular*}
	\medskip 
	\parbox{\textwidth}{\footnotesize{\textit{Notes}. The true number of staleness factors is 2. Panel A assesses the Staleness Factor Model (SFM) estimation via Root Mean Squared Error (RMSE), comparing our MLE-based estimators for the latent index ($z$) and staleness probability ($p$) against a local block (LB) method. No.Factor reports the average estimated number of factors. Panels B-D evaluate the accuracy of volatility matrix (annualized) estimation for three scenarios: (B) from efficient prices (no staleness); (C) from stale prices without correction; and (D) from stale prices with bias correction. In Panel D, estimators are denoted as follows: ($\star$) for our SFM-based MLE correction; (O) for the oracle correction using true staleness probabilities; and (LB) for the correction using probabilities from the local block method. Note that $V_t^{(p)}=V_t\circ\mathcal{P}_t$ and $\Sigma^{(p)}=\int_0^TV_t^{(p)}dt$.}}
\end{table}

The simulation results are shown in Table \ref{tab:Numerical simulation results}, providing strong support for the proposed MLE method. Panel A evaluates the finite-sample performance of the SFM estimators. Under all specifications, our MLE estimates, for both the latent index ($z$) and the staleness probability ($p$), have substantially smaller root mean squared error (RMSE) than the local-block estimates. This highlights the efficiency gains achieved by utilizing the factor structure. The estimated number of factors is also close to the true value of 2, confirming the reliability of our factor selection procedure.

Panel B-D show the impact of price staleness on volatility estimates and the effectiveness of our correction method. As shown in Panel C, the volatility estimates converge to the biased targets $V_t^{(p)}$ and $\Sigma^{(p)}$ rather than the efficient quantities $V_t$ and $\Sigma$, with a more pronounced effect for the integrated volatility matrix. Panel D shows that all correction methods mitigate this bias. The MLE-based correction outperforms the local block correction (LB) in most cases and is very close to the oracle estimator (i.e., using the true staleness probability during correction). This superior performance is consistent across link functions (Logit/Probit), sampling frequencies, and cross-sectional dimensions, particularly for the integrated volatility matrix, providing further evidence of the robustness of the proposed method.

\section{Empirical Application}\label{sec:Empirical Application}
In this section, we examine the role of staleness information in asset pricing and in estimating the volatility matrices. We sample S\&P 500 constituents at 5-minute to mitigate microstructure noise. We include the trading volume---transformed as $\log(\text{volume}+1)$---as the sole covariate. We employ a high-frequency version of the Fama–French three factors and the Carhart momentum factor, following \cite{pelger2020understanding}.\footnote{We utilize the publicly available dataset from \cite{pelger2020understanding}; https://doi.org/10.1111/jofi.12898.} Detailed data selection and cleaning procedures are documented in the supplementary materials.

\subsection{Estimation Results}\label{subsec:Estimation results}
We estimate the SFM using zero-return data from 2014 to illustrate the behavior of the staleness factors. These factors were extracted using a logit link function. We identify three distinct staleness factors, reported in Figure \ref{fig:staleness factors}. These factors exhibit markedly different daily patterns.
\begin{figure}[htb]
	\vspace{0pt}
	\includegraphics[scale=0.52]{SF_logit5.eps}
	\caption{\footnotesize Average daily staleness factors. \textit{Notes}. This figure plots the three estimated staleness factors (daily average) in 2014, based on 5-minute sampling.}\label{fig:staleness factors}
\end{figure}

Figure 15 of the supplementary material shows three representative staleness-probability trajectories exhibiting a clear comovement pattern, which ranges from a minimum of approximately 0.02 to a maximum of around 0.60. In Supplementary Figure 15 (middle panel), more than half of the stocks have staleness probabilities exceeding 0.10, indicating that staleness is pervasive in the market.

\subsection{Application in Asset Pricing}\label{subsec:Application in asset pricing}
The no-arbitrage pricing framework establishes a connection between the factors driving asset comovements and the cross-section of expected returns. In this study, we extend existing research by introducing a staleness factor to help explain the cross-sectional variation in expected excess returns. \citet{pelger2020understanding} evaluates the pricing performance of four continuous high-frequency factors against the traditional Fama-French-Carhart factors. In this section, we compare the explanatory power of the staleness factor with both sets of factors.

To effectively compare two sets of factors, we employ the generalized (canonical) correlation coefficient, following the approach of \citet{bai2006evaluating} and \citet{pelger2019large}. This measure quantifies the degree of alignment between the vector spaces spanned by two sets of factors. A coefficient of one indicates that the two factor matrices span the same subspace, while lower values reflect the highest achievable correlation between any linear combinations of the two sets. We also report canonical correlations between the staleness factors and the full panel of stock returns, providing a measure of the extent to which the staleness factors capture common variation in asset returns.
\begin{figure}[htb]
	\vspace{0pt}
	\includegraphics[scale=0.52]{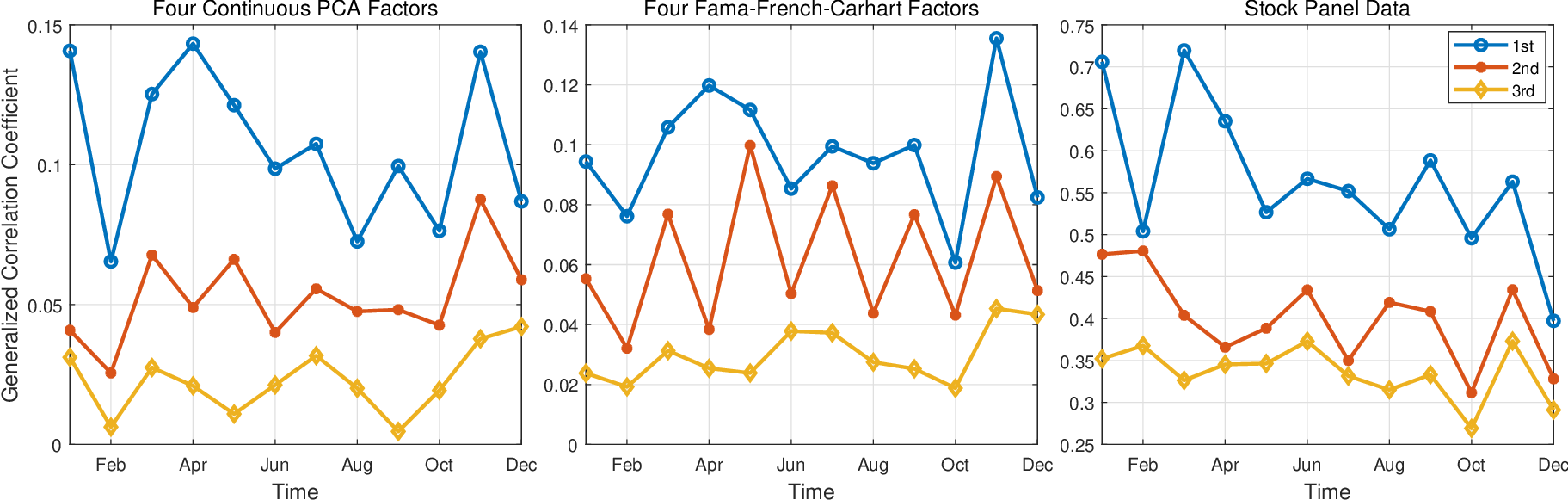}
	\caption{\footnotesize Generalized correlations between staleness factors and other factors. \textit{Notes}. The figure displays the generalized correlations of the first three staleness factors with the following: the left panel shows the four high-frequency continuous factors; the middle panel shows the Fama-French-Carhart factors; and the right panel shows the full stock-panel data. Each correlation is computed using factor estimates from a rolling one-month window throughout 2014.}\label{fig:Correlation}
\end{figure}

Figure \ref{fig:Correlation} demonstrates that the staleness factors exhibit low canonical correlations with both the high-frequency continuous factors and the Fama-French-Carhart factors, all below 0.15 throughout the sample. In contrast, the staleness factors show strong correlations with the full stock-panel data. This finding suggests that the staleness factors capture unique information inherent in the stock panel that is not reflected in the continuous or Fama-French-Carhart factor sets.

To further illustrate this point, we analyze how the proportion of variation explained by our factors evolves over time. We employ the two-stage regression framework of \cite{fama1973risk}, as extended by \cite{bollerslev2016roughing} and \cite{ait2025continuous}. Let $X_t$ denote the vector of selected factors. We conduct two comparative experiments: (i) $X_t = (FFC_t, g_t)$ versus the benchmark $X_t = FFC_t$, where $FFC_t$ represents the four Fama-French-Carhart factors, and (ii) $X_t = (CF_t, g_t)$ versus the benchmark $X_t = CF_t$, where $CF_t$ represents the four continuous factors.

\begin{figure}[htb]
	\vspace{0pt}
	\centering
	\includegraphics[scale=0.59]{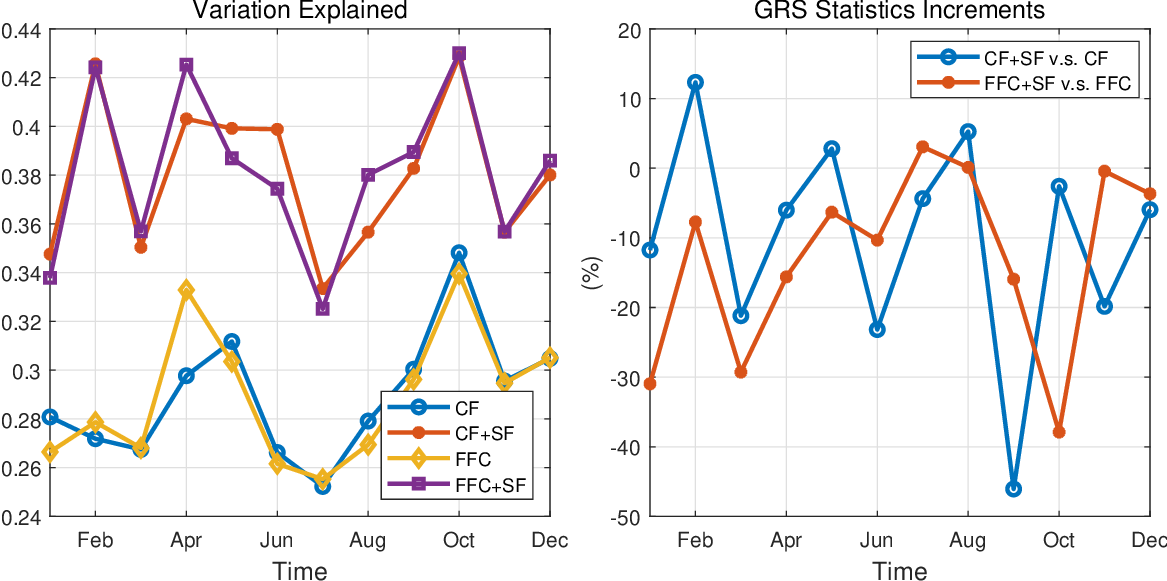}
	\caption{\footnotesize Time-varying explained variation by factor and GRS test results. \textit{Notes}. The left panel shows the percentage of continuous variation explained---computed using \cite{pelger2019large}'s method---over a rolling one-month window (21 trading days). The right panel compares the GRS statistics for the model incorporating the staleness factor with the baseline model.}\label{fig:APT}
\end{figure}

Figure \ref{fig:APT} reveals that both the four continuous factors and the Fama-French-Carhart factors explain a similar and relatively limited share of total risk. However, when the staleness factor is added to either factor set, the proportion of explained variation increases by nearly 50\%. This substantial improvement indicates that both the continuous and Fama-French-Carhart models omit critical information related to price frictions, and the staleness factor effectively captures this missing component. Additionally, the results of the Gibbons-Ross-Shanken (GRS) statistic increment (where values below zero indicate superiority over the benchmark) suggest that incorporating the staleness factor contributes to further explaining pricing errors (alpha).

\subsection{Out-of-Sample Portfolio Allocation}\label{subsec:Out-of-sample portfolio allocation}
The staleness probability can also be used to adjust the estimated volatility matrix, which may otherwise be distorted when zero returns are omitted (\citealt{kong2025data}). We assess how high-frequency, large-dimensional volatility estimates affect out-of-sample portfolio allocation by solving the constrained minimum-variance problem
\begin{align}\label{eq:portfolio allocation}
	\min_w w'\widehat{\text{cov}} w,\quad\text{s.t.}\quad w'\mathbf{1}_d=1\text{ and }\|w\|_1\leq c,
\end{align}
where $c$ is the gross-exposure bound (\citealt{fan2012vast}), and $\widehat{\text{cov}}$ denotes the (annualized) estimated covariance matrix. When $c=1$, the constraint implies no short-selling. When $c>1$, negative weights are allowed. We compare portfolios constructed using different volatility matrices---spot, integrated, and their staleness-corrected counterparts---across a range of $c$. For May 2014, we estimate $\widehat{\text{cov}}$ using April 2014 data and treat it as a plug-in forecast for the next month's covariance matrix. This design evaluates the practical benefits of correcting zero-return biases in high-frequency volatility estimation.

\begin{figure}[htb]
	\vspace{0pt}\centering
	\includegraphics[scale=0.55]{Annualized_Risk_logit5.eps}
	\includegraphics[scale=0.55]{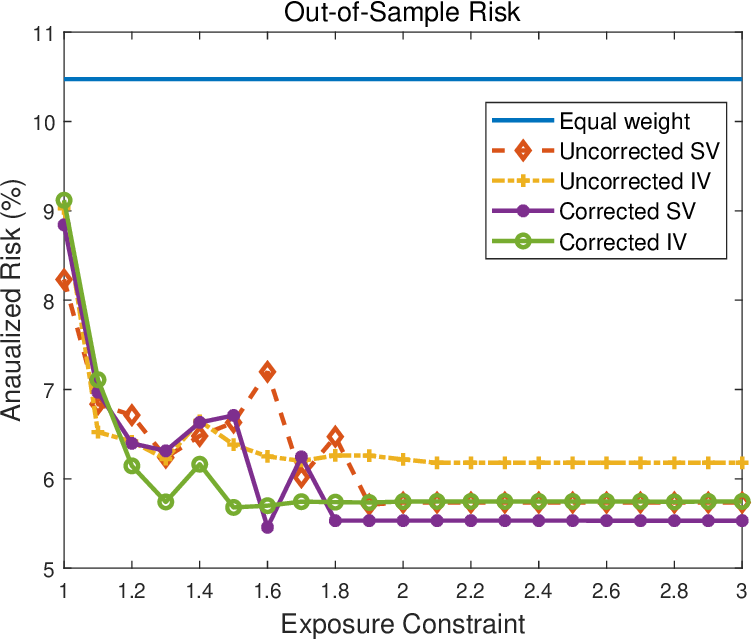}
	\caption{\footnotesize  Out-of-sample portfolio risk (left panel: 5 minute; right panel: 1 minute). \textit{Notes}. This figure compares the out-of-sample annualized volatility (for May 2014) of S\&P 500 index constituents from April 2014. The x-axis represents the exposure constraint $c$ in the optimization problem \eqref{eq:portfolio allocation}. Four volatility matrix estimators are compared: uncorrected spot volatility (Uncorrected SV), uncorrected integrated volatility (Uncorrected IV), corrected (logit type) spot volatility (Corrected SV), and corrected integrated volatility (Corrected IV). ``Equal weight'' refers to an equally weighted portfolio.}\label{fig:annualized risk}
\end{figure}

Figure \ref{fig:annualized risk} plots out-of-sample annualized risk against the gross-exposure bound $c$. For reference, we include an equal-weight portfolio---unconstrained by $c$---which exhibits a 10.5\% annualized risk. When $c=1$, the no-short-sale portfolios are poorly diversified, leading to higher out-of-sample risks. As the constraint relaxes ($c$ increases), risk declines for all estimators before leveling off.

Adjusting volatility matrices for staleness further reduces portfolio risk---especially for the integrated volatility estimator. At higher exposure levels, the staleness-corrected integrated volatility cuts risk by about 10\% compared to its uncorrected counterpart.

\section{Conclusion and Discussion}\label{sec:Conclusion}

This paper studies the cross-sectional dependence of price staleness in a general continuous-time nonlinear factor model. We introduce a novel high-frequency maximum likelihood estimation procedure and establish its asymptotic theory. We derive asymptotic bias results showing that conventional volatility matrix estimators are downward biased under price staleness, which allows us to recover and validate the latent effective price volatility matrix.

Several avenues for future research deserve exploration. First, our model assumes constant staleness factor loadings. Allowing these loadings to vary over time would be a valuable extension, though challenging because staleness manifests as binary indicators, unlike continuous value of price data. Second, we assume independence between volatility and staleness of effective prices. Exploring potential correlations between these two could yield deeper insights. 

As in many studies, we work with data at the 1-minute frequency to mitigate microstructure noise bias and to effectively filter out very high-frequency jumps in volatility estimation, price staleness, microstructure noise, and jumps may coexist and become non-negligible at ultra-high sampling frequencies (e.g., seconds or tick-by-tick). In this case, the first part of the present paper, i.e., the cross-sectional modeling of the price staleness probability, is irrelevant to the presence of price jumps and microstructure noise. This is because the price staleness defined in the literature relates only to the occurrence of flat prices (or equivalently zero returns) which are completely observed and not attributed to the jumps and noise. But the second part, estimating the volatility matrix is definitely biased due to the presence of the jumps and noise, besides the price staleness revealed in the present paper. In the ultra-high frequency setting, following standard practice in the literature, for each asset, one could remove the noise by locally smoothing the data (e.g., pre-averaging) and then truncate the jumps (large smoothed increments), which produces a nearly noise-free, jump-robust, and approximately linearly transformed increments of the continuous martingale component. With the approximate transformed continuous increments, a possible way is to implement the POET estimator of a large volatility matrix and then debias the estimates of the co-volatilities due to the price staleness under model \eqref{eq:observable price processes}. The theoretical perspective of this three-step (pre-averaging + jump truncation + staleness correction) sequential approach needs more lengthy and delicate mathematical analysis. Conceptually, the price staleness is a kind of market microstructure that is different from the stylized additive microstructure noise on top of the efficient prices. So whether there is a uniform method that can remove the bias due to the additive noise and the price staleness from estimating the co-volatilities of the continuous components of assets is an interesting future research topic.






\bibliographystyle{Chicago}
\spacingset{1.4}
\bibliography{ref}
\end{document}